\documentclass[11pt]{amsart}
\usepackage{amsmath}
\usepackage{amssymb}
\usepackage{amsfonts}
\usepackage{geometry}
\usepackage{leqno}
\usepackage{mathtools}
\usepackage{amsthm}
\usepackage[shortlabels]{enumitem}
\usepackage{xcolor}
\usepackage[all]{xy}
\usepackage{nicefrac}
\definecolor{cobalt}{RGB}{61,89,171}
\usepackage[colorlinks,citecolor=cobalt,linkcolor=cobalt,urlcolor=cobalt,pdfpagemode=UseNone,backref = page]{hyperref}

\newcommand{\hook}{\lrcorner \,}

\newcommand\mfa{{\mathfrak a}}
\newcommand\frg{{\mathfrak g}}
\newcommand\mfg{{\mathfrak g}}
\newcommand\frh{{\mathfrak h}}
\newcommand\mfh{{\mathfrak h}}
\newcommand\frj{{\mathfrak j}}
\newcommand\mfj{{\mathfrak j}}
\newcommand\frl{{\mathfrak l}}
\newcommand\mfr{{\mathfrak r}}
\newcommand\mfu{{\mathfrak u}}

\newcommand{\mfup}{\mathfrak{u}^{\perp_\omega}}
\newcommand\mfz{{\mathfrak z}}
\newcommand\hrho{{\hat{\rho}}}
\newcommand{\bJ}{\mathbf{J}}
\newcommand{\bO}{\mathbf{\Omega}}
\newcommand{\bC}{\mathbb{C}}
\newcommand{\bN}{\mathbb{N}}
\newcommand{\bQ}{\mathbb{Q}}
\newcommand{\bR}{\mathbb{R}}
\newcommand{\diag}{{\mathrm{diag}}}
\newcommand{\ad}{\operatorname{ad}}
\newcommand{\id}{\operatorname{id}}
\newcommand{\End}{\operatorname{End}}
\newcommand{\GL}{\operatorname{GL}}
\newcommand{\im}{\operatorname{Im}}
\newcommand{\re}{\operatorname{Re}}

\newcommand{\spa}[1]{\mathrm{span}(#1)}
\newcommand{\Ann}[1]{#1^0}

\newtheorem{theorem}{Theorem}[section]
\newtheorem*{theorem*}{Theorem}
\newtheorem{lemma}[theorem]{Lemma}
\newtheorem{corollary}[theorem]{Corollary}
\newtheorem{proposition}[theorem]{Proposition}
\newtheorem{claim}[theorem]{Claim}
\newtheorem{example}[theorem]{Example}

\newtheorem{definition}[theorem]{Definition}

\newtheorem{notation}[theorem]{Notation}
\theoremstyle{remark}
\newtheorem{remark}[theorem]{Remark}

\title{Complex symplectic Lie algebras with large Abelian subalgebras}

\author[G. Bazzoni]{Giovanni Bazzoni}
\address{Dipartimento di Scienza ed Alta Tecnologia, Universit\`a degli Studi dell'Insubria, Via Valleggio 11, 22100, Como, Italy}
\email{giovanni.bazzoni@uninsubria.it}

\author[M. Freibert]{Marco Freibert}
\address{Mathematisches Seminar\\
Christian-Albrechts-Universit\"at zu Kiel\\
Heinrich-Hecht-Platz 6\\
D-24118 Kiel\\
Germany}

\email{freibert@math.uni-kiel.de}

\author[A. Latorre]{Adela Latorre}
\address{Departamento de Matem\'atica Aplicada\\
Universidad Polit\'ecnica de Madrid\\
Avda. Juan de Herrera 4\\
28040 Madrid, Spain}
\email{adela.latorre@upm.es}

\author[N. Tardini]{Nicoletta Tardini}
\address{Dipartimento di Scienze Matematiche, Fisiche e Informatiche\\
Unit\`a di Matematica e Informatica\\
Universit\`a degli Studi di Parma\\
Parco Area delle Scienze 53/A\\
43124 Parma, Italy}
\email{nicoletta.tardini@unipr.it}

\begin{document}


\maketitle

\begin{abstract}
We present two constructions of complex symplectic structures on Lie algebras with large Abelian ideals. In particular, we completely classify complex symplectic structures on almost Abelian Lie algebras.
By considering compact quotients of their corresponding connected, simply connected Lie groups we obtain many examples of complex symplectic manifolds which do not carry (hyper)k\"ahler metrics. We also produce examples of compact complex symplectic manifolds endowed with a fibration whose fibers are Lagrangian tori.
\end{abstract}

\section{Introduction}

A complex symplectic structure on a complex manifold $X$ is a holomorphic 2 form $\sigma$ which is closed and non-degenerate. Hence the complex dimension of $X$ is even, say $2n$, and the canonical bundle of $X$ is trivialized by $\sigma^n$. By the celebrated Newlander-Nirenberg theorem, one can equivalently think of $X$ as a pair $(M,J)$ where $M$ is a smooth manifold and $J$ is an integrable almost complex structure, that is, an endomorphism of $TM$ with $J^2=-Id$ and $N_J=0$, where $N_J$ is the Nijenhuis tensor. In this language, a complex symplectic structure is a $(2,0)$-form $\sigma\in\Omega^{2,0}(M,J)$ which is closed and non-degenerate. This turns out to be equivalent to a symmetry condition on $\omega=\sigma+\bar\sigma$, namely $\omega(JX,Y)=\omega(X,JY)$, for all vector fields $X,Y$, see Lemma \ref{le:equivalencecomplexsymplectic}. This compatibility condition between $\omega$ and $J$ is antithetical to the one defining a K\"ahler metric. Similarly to the {\em real} symplectic situation, complex symplectic structures exist on the holomorphic cotangent bundle of any complex manifold, and on coadjoint orbits of complex Lie groups \cite{Boalch}; such examples, however, are not compact. Hyperk\"ahler geometry is a source of examples of complex symplectic manifolds: if $(g,I_1,I_2,I_3)$ is a hyperk\"ahler structure on a smooth manifold $M$, with symplectic forms $\omega_j=g(\cdot,I_j\cdot)$, $j=1,2,3$, then $\sigma=\omega_2+i\omega_3$ is a complex symplectic structure on $(M,I_1)$. In fact, by a result of Beauville \cite{Beauville}, a compact complex manifold of K\"ahler type admits a hyperk\"ahler structure if and only if it admits a complex symplectic structure.  Hyperk\"ahler manifolds are important both from the Riemannian and the complex point of view, see for instance \cite{GHJ,Huy,Joyce,Sal}; however compact examples of these manifolds do not abound. In \cite{Guan1,Guan2}, Guan constructed compact complex symplectic non (hyper)k\"ahlerian manifolds. His constructions use nilmanifolds, compact quotients of a connected, simply connected nilpotent Lie group by a lattice, in a crucial way. Natural complex symplectic structures also exist on Hitchin's hypersymplectic manifolds, which are the closest relatives of hyperk\"ahler manifolds in the context of neutral signature, see \cite{Dancer-Swann,Hitchin}.

An open question in real symplectic geometry concerns the existence of symplectic structures on compact manifolds; therefore, being able to construct examples is essential. Also, a great deal of research over the last years has addressed the question of the existence of compact symplectic manifolds with no K\"ahler metrics - see for instance \cite{BFM,CFM,FM,Tralle-Oprea}. Two natural questions arise at this point:
\begin{itemize}
\item How to construct compact complex symplectic manifolds?
\item In case one is constructed, how to exclude that it carries a - say - hyperk\"ahler metric?
\end{itemize}
In \cite{BFLM}  the first three authors addressed both questions studying the existence of complex symplectic structures on (non toral) nilmanifolds, 
which are known to admit no (hyper)K\"ahler structure, see \cite{Benson-Gordon, Hasegawa}. The structures considered there are left-invariant, hence defined at the Lie algebra level. In particular, complex symplectic oxidation is introduced in order to provide a method to construct $(4n+4)$-dimensional complex symplectic Lie algebras (not necessarily nilpotent) 
starting from a $4n$-dimensional one. The construction is based on the existence of a $2$-dimensional central $J$-invariant ideal. In this paper, which is a kind of follow-up to \cite{BFLM}, we tackle again both questions
and describe two new constructions of complex symplectic structures on solvable Lie algebras. The leitmotiv behind the title is that the kind of Lie algebras we consider have ``large'' Abelian ideals: of codimension one in the first case, and of half-dimension in the second one.
These structures induce left-invariant complex symplectic structures on the corresponding connected, simply connected solvable Lie groups and on their solvmanifolds. We thus provide systematic ways of constructing examples of compact complex symplectic (solv)manifolds, and show that they abound. As we recalled, no nilmanifolds (except tori) admit K\"ahler metrics, and solvmanifolds with K\"ahler metrics are completely understood, see \cite{Hasegawa2}. Thus we can exclude that our examples admit K\"ahler, hence hyperk\"ahler, metrics. We also study solvable Lie algebras with Abelian complex symplectic structures, those whose complex structure is Abelian. This also fits in the title, since an Abelian complex structure on a Lie algebra $\frg$ is given by an Abelian subalgebra $\frg^{(1,0)}\subset\frg\otimes\bC$. 

Here is an overview of the content and of the main results in this paper:
\begin{itemize}
\item In Section \ref{sec:Preliminaries} we recall some preliminaries and fix the notation.
\item In Section \ref{sec:almost_Abelian} we obtain all almost Abelian Lie algebras with a complex symplectic structure, see Theorem \ref{th:complexsymplectic}. To attain this result, we first recall that the integrability of an almost complex structure~$J$ on an almost Abelian Lie algebra $\frg$ can be characterized by the (matrix) structure of $f\coloneqq\ad_X|_{\mfu}$, being~$\mfu$ the codimension one Abelian ideal of $\frg$ and $X\in\frg\setminus\mfu$. We then provide a result in similar terms to characterize when a non-degenerate two-form~$\omega$ on~$\frg$ is closed. Note that, at this point of the paper (Section~\ref{subsec:complex-and-symplectic}), the complex and the symplectic structures are non-related. Hence, the next step is to make them interact (Section~\ref{subsec:complexsymplecticalmostAbelian}). We show that the two previous descriptions for $f$ can be combined when we impose that $(J,\omega)$ is a complex symplectic structure. We make use of these results to show that, in general, an almost Abelian Lie algebra can admit two non-isomorphic complex symplectic structures, see Example \ref{ex:nonuniqueness}.
\item Section \ref{sec:classification} contains the classification of almost Abelian Lie algebras admitting complex symplectic structures, characterized in terms of the Jordan blocks of $f$, see Theorem \ref{th:almostAbelianLAclass}. 
More precisely, we first observe that any complex symplectic almost abelian Lie algebra is given by $(\mathbb R^{4n-1}\rtimes_f\mathbb R,J_0,\omega_0)$, where $f$ is determined by Theorem~\ref{th:complexsymplectic} and $(J_0,\omega_0)$ takes a canonical form. Then, the equivalence classes of $f$ that preserve the canonical complex symplectic structure are studied. This allows to provide the set of matrices for $f$ that give rise to non-equivalent complex symplectic almost Abelian Lie algebras (Proposition~\ref{pro:classsimplification}), from where one finally derives the theorem. In Corollary \ref{co:uniqueness} we obtain conditions under which an almost Abelian Lie algebra admits a unique complex symplectic structure. In \ref{subsec:solvmanifold} we also construct almost Abelian solvmanifolds with a unique left-invariant complex symplectic structure, and which carry no K\"ahler metric.
\item In Section \ref{sec:cotangent} we describe the cotangent extension, which starts with a Lie algebra $\frg$ endowed with a complex structure $J$ and gives back a complex symplectic structure on $\frh=\frg^*\oplus\frg$, endowed with a Lie algebra structure such that $\frg^*\subset\frh$ is a complex, Abelian, Lagrangian ideal. By Theorem \ref{theo:cot_extension}, this construction completely characterizes such complex symplectic Lie algebras. A geometric version of this result yields complex symplectic manifolds endowed with a Lagrangian fibration, see \ref{subsec:Lag_fib}.
\item Section \ref{sec:Abelian} deals with complex symplectic structures on solvable Lie algebras with non-trivial center under the assumption that the complex structure is Abelian. In Proposition \ref{pro:conditionsJAbelian} we obtain conditions under which the reduction-oxidation procedure of \cite{BFLM} produces such Lie algebras. Using this conditions in Theorem \ref{thm:nilp_Abelian} we prove that nilpotent Lie algebras with Abelian complex symplectic structures can be obtained by iterated oxidations of the trivial Lie algebra. Moreover, we show that, for nilpotent Lie algebras with Abelian complex symplectic structures, all the possible nilpotency steps are admissible. More precisely, in Proposition \ref{pro:nilpotencystepAbelian} we show that for a fixed $n\in \bN$, then there exists a $4n$-dimensional nilpotent complex symplectic Lie algebra of step length $m$, for any $m\in \{1,\ldots,2n\}$.
\end{itemize}

\medskip 
\noindent {\bf Acknowledgements.} The first author is partially supported by GNSAGA and by grants PID2020-118452GB-I00 and  PID2021-126124NB-I00 (MCIN/AEI/10.13039/501100011033). The third author is partially supported by grant PID2020-115652GB-I00, funded by MCIN/AEI/10.13039/501100011033. The fourth author has financially been supported by the Programme ``FIL-Quota Incentivante'' of University of Parma and co-sponsored by Fondazione Cariparma and by GNSAGA of INdAM. We are indebted to the referee for her/his comments, which helped us improving the presentation of the paper.

\section{Preliminaries}\label{sec:Preliminaries}

In this section we present the basic concepts that will be used in the paper, and fix the notation. Let $M$ be a $2m$-dimensional differentiable manifold, and let $\mathfrak X(M)$ denote the space of smooth vector fields on $M$. A \emph{complex structure} $J$ on $M$ is an endomorphism $J:\mathfrak X(M)\to\mathfrak X(M)$ satisfying $J^2=-Id$ and $N_J(X,Y)=0$ for every $X,Y\in\mathfrak X(M)$, where $N_J$ is the Nijenhuis tensor
$$N_J(X,Y):=[X,Y]+J[JX,Y]+J[X,JY]-[JX,JY].$$
The pair $(M,J)$ is then a \emph{complex manifold} of complex dimension $m$, as a consequence of the celebrated Newlander-Nirenberg theorem. $J$ can be equivalently defined on $\Omega^1(M)$, the space of smooth $1$-forms on~$M$, by taking
\begin{equation}\label{J-on-forms}
	(J^*\alpha)(X)=\alpha(JX),
\end{equation}
for every $\alpha\in\Omega^1(M)$ and $X\in\mathfrak X(M)$. Extending $J^*$ by $\mathbb C$-linearity to $\Omega^1_{\mathbb C}(M)=\Omega^1(M)\otimes\mathbb C$, one obtains a decomposition
$$\Omega^1_{\mathbb C}(M)=\Omega^{1,0}(M,J)\oplus\Omega^{0,1}(M,J),$$
where 
\begin{equation*}
\begin{split}
\Omega^{1,0}(M,J)&=\left\{\alpha\in\Omega^1_{\mathbb C}(M)\mid J^*\alpha=i\,\alpha \right\},\\
\Omega^{0,1}(M,J)&=\{\alpha\in\Omega^1_{\mathbb C}(M)\mid J^*\alpha=-i\,\alpha \}
	= \overline{\Omega^{1,0}(M,J)}.
\end{split}
\end{equation*}
As a consequence, the space of complex smooth $k$-forms splits as
$$\Omega^k_{\mathbb C}(M)=\bigoplus_{p+q=k}\Omega^{p,q}(M,J),$$
where $\Omega^{p,q}(M,J)$ denotes the space of $(p,q)$-forms on $(M,J)$.

A \emph{complex symplectic}, or \emph{holomorphic symplectic}, structure on $M$ is a pair $(J,\omega_{\mathbb C})$ given by a complex structure $J$ on $M$ and a non-degenerate $\omega_{\mathbb C}\in\Omega^{2,0}(M,J)$ satisfying $d\omega_{\mathbb C}=0$. Note that the existence of such a $(J,\omega_{\mathbb C})$ forces the complex dimension of $M$ to be even, namely, $m=2n$. Moreover, the $(2,0)$-form $\omega_{\mathbb C}$ can be replaced by a (real) symplectic form satisfying some additional condition, as shown in \cite[Lemma 3.2]{BFLM}:

\begin{lemma}\label{le:equivalencecomplexsymplectic}
Let $M$ be a differentiable manifold. The set of all complex symplectic structures $(J,\omega_\bC)$ on $M$ is bijective to the set of all pairs $(J,\omega)$ consisting of a complex structure $J$ and a symplectic 2-form $\omega$ such that $J$ is symmetric with respect to $\omega$, i.e.~$\omega(JX,Y)=\omega(X,JY)$ for all $X,Y\in \mathfrak{X}(M)$. The bijection is given by $(J,\omega_\bC)\mapsto (J,\Re(\omega_\bC))$ with inverse $(J,\omega)\mapsto (J,\omega-i \omega(J\cdot,\cdot))$.
\end{lemma}
In this paper we consider the case in which $M$ is a \emph{solvmanifold}, namely, a quotient $M=\Gamma\backslash G$ of a real, connected, simply connected solvable Lie group $G$ by a lattice~$\Gamma$ (i.e.~a discrete and cocompact subgroup). {\em Nilmanifolds} are special solvmanifolds, those for which $G$ is nilpotent. When $G$ is nilpotent, the existence of lattice $\Gamma\subset G$ is equivalent to the Lie algebra $\frg$ of $G$ having a {\em rational structure}, that is, a rational subalgebra $\mfg_\bQ\subset \mfg$ such that $\mfg=\mfg_\bQ\otimes\bR$ \cite{malcev}. There is however no sufficient condition for the existence of lattices in solvable Lie groups: their construction is art. {\em Almost Abelian} solvable Lie groups form a little but important exception to this difficulty, see \cite{bock} and Subsection \ref{subsec:solvmanifold} below.

Next, recall that a Lie algebra can also be described using $d\colon\frg^*\to\Lambda^2\frg^*$, the dual map to the Lie bracket, defined by $(d\alpha)(X,Y)=-\alpha([X,Y])$. The Jacobi identity on $\frg$ is then equivalent to $d^2=0$, where $d$ is extended to $\Lambda^*\frg^*$ as a graded derivation. Suppose $\{e_1,\ldots,e_n\}$ is a basis of $\frg$, and that $[e_i,e_j]=\sum c^k_{ij}e_k$. If $\{e^1,\ldots,e^n\}$ is the dual basis, then $de^k=-\sum c^k_{ij}e^i\wedge e^j$. We shall use the so-called {\em Salamon's notation} for denoting Lie algebras. For instance, $\frg=(0^3,12,15,-16)$ means that the 6-dimensional Lie algebra $\frg$ has a basis $\{e_1,\ldots,e_6\}$ with dual basis $\{e^1,\ldots,e^6\}$ such that $de^i=0$, $i=1,2,3$, $de^4=e^{12}\coloneqq e^1\wedge e^2$, $de^5=e^{15}$ and $de^6=-e^{16}$.


A complex structure on a solvmanifold $M=\Gamma\backslash G$ is {\em invariant} if it comes from a left-invariant complex structure on $G$, that is, from a complex structure on $\frg=\textrm{Lie}(G)$; this is just an endomorphism $J\colon\frg\to\frg$ with $J^2=-Id$ and $N_J=0$, and is equivalently described by a subspace $\frg^{*(1,0)}\subset\frg^*\otimes\bC$, precisely as it happens on smooth manifolds. There are distinguished classes of complex structures on Lie algebras.

\begin{definition}\label{def-types-J}
Let $\frg$ be a real Lie algebra endowed with an almost complex structure $J$.
\begin{itemize}
\item $J$ is {\em Abelian} if $[X,Y]=[JX,JY]$ for all $X,Y\in\frg$; equivalently, $J$ is Abelian if $d(\frg^{*(1,0)})\subset\frg^{*(1,1)}$ if and only if $\frg^{(1,0)}\subset\frg\otimes\bC$ is an Abelian subalgebra;
\item $J$ is {\em parallelizable} if $J[X,Y]=[JX,Y]$ for all $X,Y\in\frg$; equivalently, $J$ is parallelizable if $d(\frg^{*(1,0)})\subset\frg^{*(2,0)}$. In this case, $G$ has the structure of a complex Lie group.
\end{itemize}
Notice that in both cases $J$ is automatically integrable. 
\end{definition}

In this paper we shall work with (solvable) Lie algebras, and then obtain geometric examples by constructing lattices in the corresponding Lie groups.

\begin{definition}
Let $\frg$ be a real Lie algebra. A {\em complex symplectic structure} on $\frg$ consists of an integrable almost complex structure $J$ and a symplectic form $\omega$ such that $J$ is symmetric with respect to $\omega$. Thus $N_J=0$, $\omega$ is non-degenerate, $d\omega=0$, and $\omega(JX,Y)=\omega(X,JY)$ for every $X,Y\in\frg$. Then, $(\mathfrak g, J, \omega)$ is called a \emph{complex symplectic Lie algebra}.
\end{definition}

Taking into account the notions of equivalence between real Lie algebras, complex structures and symplectic structures, it is natural to introduce the following:

\begin{definition}\label{def:equivalence-CS}
	Two complex symplectic Lie algebras $(\mathfrak g, J, \omega)$ and
	$(\mathfrak g', J', \omega')$ are said to be \emph{equivalent} if there exists an isomorphism of Lie algebras $\phi:\mathfrak g\to\mathfrak g'$
	satisfying $\phi\circ J=J'\circ \phi$ and 
	$\phi^*\omega'=\omega$.
\end{definition}

\section{Complex symplectic structures on almost Abelian Lie algebras}\label{sec:almost_Abelian}

In this section we construct complex symplectic structures on real almost Abelian Lie algebras. We start by recalling some generalities and proving some results concerning the existence of complex structures and symplectic structures on this type of Lie algebras. Then, we focus on the complex symplectic case.

\begin{definition}
A (solvable) real Lie algebra $\frg$ is {\em almost Abelian} if it has a codimension one Abelian ideal. A Lie group $G$ is almost Abelian if its Lie algebra is.
\end{definition}



Recall the following result in \cite{avetisyan}:

\begin{proposition}\label{pro:uniquenesscodim1ideal}
Let $\mathfrak h$ be an almost Abelian Lie algebra over the field $\mathbb F\in\{\mathbb R,\mathbb C\}$. If $\mathfrak h$ has more than one codimension one Abelian ideals, then $\mathfrak h\cong \mathbb H_{\mathbb F}\oplus\mathbb F^k$ for some $k\in\mathbb N_0$, where $\mathbb H_{\mathbb F}$ is the Heisenberg algebra
$$\mathbb H_{\mathbb F}=\left\{
	\left.\begin{pmatrix} 0 & 0 & 0 \\ t & 0 & p\\ q & 0 & 0 \end{pmatrix}
	\quad \right\vert \quad
	(p,t,q)\in\mathbb F^3
\right\}.$$
\end{proposition}


Note that the Lie algebras in the result above are nilpotent. Indeed,
there are two well-known complex nilmanifolds whose associated Lie algebras fit into the previous statement.

If we take $\mathbb F=\mathbb R$, then $\mathbb H_{\mathbb R}$ is usually denoted by $\mathfrak h_3$. $\mathfrak r\mathfrak h_3\coloneqq\mathfrak h_3\oplus\mathbb R$ is the Lie algebra underlying the Kodaira-Thurston nilmanifold, which was 
the first example of a compact manifold with both symplectic and complex structures but no K\"ahler metric \cite{thurston}. It has been shown in \cite{BGGL} that there is only one complex symplectic structure on $\mathfrak r\mathfrak h_3$ up to equivalence. 

$\mathbb H_{\mathbb C}$ is the complex Lie algebra underlying the Iwasawa manifold. For dimension reasons, it does not admit a complex symplectic structure; however $\mathfrak h\coloneqq\mathbb H_{\mathbb C}\oplus\mathbb C$ is determined by the complex structure equations
\begin{equation}\label{iwasawa}
d\omega^1=d\omega^2=0\,,\qquad d\omega^3=\omega^{12}\,, \qquad d\omega^4=0\,.
\end{equation}
It is easy to see that any closed non-degenerate complex $2$-form $\omega_{\mathbb C}$ on $\mathfrak h$ is given by
$$\omega_{\mathbb C}=a_{12}\omega^{12}+a_{13}\omega^{13}+a_{14}\omega^{14}
+a_{23}\omega^{23}+a_{24}\omega^{24},$$
with $a_{ij}\in\mathbb C$ and $a_{13}a_{24}-a_{14}a_{23}\neq 0$. Moreover, one can define a new (complex) basis for $\mathfrak h$ as follows
\begin{equation*}
\begin{split}
\tau^1 &= \frac{1}{(a_{13}a_{24}-a_{14}a_{23})^{\nicefrac 13}}\,(a_{13}\omega^1+a_{23}\omega^2),\\
\tau^2 &= \frac{1}{(a_{13}a_{24}-a_{14}a_{23})^{\nicefrac 13}}\,(a_{14}\omega^1+a_{24}\omega^2),\\
\tau^3 &= (a_{13}a_{24}-a_{14}a_{23})^{\nicefrac 13}\,\omega^3,\\
\tau^4 &= (a_{13}a_{24}-a_{14}a_{23})^{\nicefrac 13}\,\omega^4
	-\frac{a_{12}}{(a_{13}a_{24}-a_{14}a_{23})^{\nicefrac 23}}\,(a_{13}\omega^1+a_{23}\omega^2).
\end{split}
\end{equation*}
A direct calculation shows that the complex structure equations for $\mathfrak h$ in terms of the new basis $\{\tau^k\}_{k=1}^4$ still follow \eqref{iwasawa}, but now 
$$\omega_{\mathbb C}=\tau^{13}+\tau^{24}.$$
Consequently, there is just one complex symplectic structure on $\mathfrak h$ up to equivalence.
\begin{remark}
In the notation used in~\cite{ugarte}, the real $6$-dimensional Lie algebra underlying $\mathbb H_{\mathbb C}$ is $\mathfrak h_5=(0,0,0,0,13-24,14+23)$, which is not almost Abelian. According to \cite{COUV}, there are more complex structures on $\mathfrak h_5$, but none of them is parallelizable (see Definition~\ref{def-types-J}).
\end{remark}

In the next sections, we will only focus on real almost Abelian Lie algebras. In particular, we will address the existence and the uniqueness of complex symplectic structures on real almost Abelian Lie algebras whose codimension one Abelian ideal is unique.

\subsection{Complex structures and symplectic structures on almost Abelian Lie algebras}\label{subsec:complex-and-symplectic}
In this subsection, we characterise the existence of 
two different geometric structures on real almost Abelian Lie algebras:
complex structures and symplectic structures. Notice that, in this part of the paper, the two structures are treated independently. Moreover, here we do not need to require the Abelian ideal to be unique. These results will then be put together in Subsection \ref{subsec:complexsymplecticalmostAbelian} to give a characterisation of complex symplectic structures on almost Abelian Lie algebras.

The case of complex structures on almost Abelian Lie algebras has been treated individually several times in the literature. We refer to \cite{LR}, which seems to be the first characterization of complex structures on almost Abelian Lie algebras.
\begin{theorem}\label{th:complex}
Let $\mfg$ be a $2n$-dimensional almost Abelian Lie algebra with codimension one Abelian ideal $\mfu$ and let $J$ be an almost complex structure on $\mfg$. Choose $X\in \mfg\setminus \mfu$ such that $JX\in \mfu$ and set $\mfu_J\coloneqq\mfu\cap J\mfu$ and $f\coloneqq\ad_X|_{\mfu}\in \End(\mfu)$. Then $J$ is integrable if and only if there are $f_0\in \mathfrak{gl}(\mfu_J,J)\coloneqq\left\{h\in \End(\mfu_J)\mid [h,J]=0\right\}\cong \mathfrak{gl}(n-1,\bC)$, $v\in \mfu_J\cong \bR^{2n-2}$ and $a\in \bR$ such that
\begin{equation*}
f=\begin{pmatrix}
      f_0 & v \\
      \mathbf 0 & a
   \end{pmatrix}
\end{equation*}
with respect to the splitting $\mfu=\mfu_J\oplus \spa{JX}$.
\end{theorem}

Next, we investigate the existence of symplectic structures on almost Abelian Lie algebras. For this, we first need the following observations.
\begin{remark}
Let $\mathfrak g$ be an almost Abelian Lie algebra of dimension $2n$ with codimension one Abelian ideal $\mathfrak u$. Let $\omega$ be a non-degenerate $2$-form on $\mathfrak g$ and denote by $\mfu^{\perp_\omega}$ the symplectic orthogonal of $\mfu$, 
\[
\mfup\coloneqq\left\{ Y\in\mfg \mid \omega(Y,U)=0, \ \forall\,U\in\mfu \right\}\,.
\]
Since $\text{dim}\,\mfu +\text{dim}\,\mfu^{\perp_\omega}=\text{dim}\,\mfg$, $\mfup$ is one-dimensional, hence isotropic. Therefore, $\mfu=\big(\mfup\big)^{\perp_\omega}$ is coisotropic and thus $\mfup\subset\mfu$. By \cite[Lemma 2.2]{McDS}, $\mfu/\mfup$ has a symplectic structure induced by $\omega$. Thus, any complement $\mfu'$ of $\mfup$ in $\mfu$ inherits a symplectic structure $\omega'$, via the isomorphism $\mfu/\mfup\cong\mfu'$. We can therefore write
\begin{equation}\label{eq:775}
\mfg=\mfu'\oplus(\mfu')^{\perp_\omega} \text{ \ and \ } (\mfu')^{\perp_\omega}=\mfup\oplus \spa X\,,
\end{equation}
for some $0\neq X\in\mfg\setminus\mfu$. Note also that $\omega(X,Y)\neq 0$ for every $0\neq Y\in\mfup$, as $\omega$ is non-degenerate.
\end{remark}

The remark allows us to state the following result:

\begin{theorem}\label{th:symplectic}
	Let $\mfg$ be a $2n$-dimensional almost Abelian Lie algebra with codimension one Abelian ideal $\mfu$ and let $\omega$ be a non-degenerate two-form on $\mfg$. 
	Consider some complement~$\mfu'$ of $\mfup$
	in $\mfu$ and set $\omega'\coloneqq\omega|_{\mfu'}\in \Lambda^2 \mfu'$. Moreover, choose $X\in \mfg\setminus \mfu$ and denote $f\coloneqq\ad_X|_{\mfu}$. Then $\omega$ is a symplectic form on $\mfg$ if and only if there exist $f'\in \mathfrak{sp}(\mfu',\omega')\coloneqq\left\{h\in \mathrm{End}(\mfu')\mid h.\omega'=0\right\}\cong \mathfrak{sp}(2n-2,\bR)$, $\alpha\in \mathrm{Hom}(\mfu',\mfup)\cong \bR^{2n-2}$ and $a'\in \bR$ such that
	\begin{equation*}
	f=\begin{pmatrix}
	 f' & \mathbf 0 \\
	 \alpha & a'
	\end{pmatrix}
	\end{equation*}
	with respect to the splitting $\mfu=\mfu'\oplus \mfup$.
\end{theorem}
\begin{proof}
First of all, note that we need to check when
\begin{equation*}
0=d\omega(Z_1,Z_2,Z_3)=-\sum_{cycl} \omega\big([Z_i,Z_{i+1}],Z_{i+2}\big)
\end{equation*}
holds for any $Z_1,Z_2,Z_3\in \mfg$. This is trivially fulfilled if $Z_1,Z_2,Z_3\in \mfu$. So we only need to check $d\omega(X,Y,Z)=0$ and $d\omega(X,Z_1,Z_2)=0$ for all $Y\in \mfup$ and all $Z,Z_1,Z_2\in \mfu'$. We first have
\begin{equation*}
d\omega(X,Y,Z)=-\omega(f(Y),Z)-\omega(Y,f(Z))=-\omega(f(Y),Z)
\end{equation*}
due to $Y\in \mfup$ and $f(Z)\in \mfu$. Hence,
$d\omega(X,Y,Z)=0$ for all $Y\in \mfup$ and $Z\in \mfu'$ is equivalent to $f(Y)\in \mfup$ for all $Y \in \mfup$, i.e. to $f(Y)=a'Y$ for some $a'\in \bR$. Finally, write $f|_{\mfu'}=f'+\alpha$ for $f'\in \End(\mfu')$ and $\alpha\in \mathrm{Hom}(\mfu',\mfup)$. 
Note that we have $\omega(\alpha(\mfu'),\mfu')=0$ and so
\begin{equation*}
\begin{split}
d\omega(X,Z_1,Z_2)&=-\omega(f(Z_1),Z_2)-\omega(Z_1,f(Z_2))=-\omega'(f'(Z_1),Z_2)-\omega'(Z_1,f'(Z_2))\\
&=(f'.\omega')(Z_1,Z_2)
\end{split}
\end{equation*}
for all $Z_1,Z_2\in \mfu'$. Hence $d\omega(X,Z_1,Z_2)=0$ for all $Z_1,Z_2\in \mfu'$ if and only if $f'\in \mathfrak{sp}(\mfu',\omega')$.
\end{proof}

\subsection{Almost Abelian Lie algebras admitting a complex symplectic structure}\label{subsec:complexsymplecticalmostAbelian}

We combine now Theorem \ref{th:complex} and Theorem \ref{th:symplectic} in order to identify almost Abelian Lie algebras admitting complex symplectic structures. 
In particular, note that if an almost Abelian Lie algebra $\mathfrak g$ has both a complex structure $J$ and a symplectic structure $\omega$, then the previous theorems provide two different splittings for $\frg$, namely, $\frg=\mfu_J\oplus\spa{JX}\oplus\spa X$ and $\frg=\mfu'\oplus\mfu^{\perp_\omega}\oplus\spa X$, being $X$ an element of $\frg$ which is not in its codimension one Abelian ideal $\mfu$. We here show that these two splittings can be combined to find a finer one, as long as the pair $(J,\omega)$ gives rise to a complex symplectic structure.

Let $\mfg$ be a $4n$-dimensional almost Abelian Lie algebra with codimension one Abelian ideal $\mfu$. Let $(J,\omega)$ be a \emph{$\mathrm{Sp}(2n,\bC)$-structure} on $\mfg$, i.e. $\omega$ is a non-degenerate two-form on $\mfg$ and $J$ is an almost complex structure on $\mfg$ which is symmetric with respect to $\omega$. Consider again the subspaces $\mfu_J\coloneqq\mfu\cap J\mfu$ and $\mfup$ of $\mfu$. 
\begin{claim}
$\mfup\subset \mfu_J$.
\end{claim}

If not, there exists $0\neq Y\in \mfup$ such that $JY\notin \mfu$, and  $\mfg=\mfu\oplus \spa{JY}$. Since $Y\in \mfu^{\perp\omega}$ and $\omega(Y,JY)=0$ due to the symmetry of $J$, this would imply that $\omega$ is degenerate. As a consequence, one also has $J\mfup\subset \mfu_J$.
\begin{flushright}
\qedsymbol{}
\end{flushright}

\medskip
Next, choose a complement $\mfu'$ of $\mfup$ in $\mfu$ containing $J\mfup$ and set $\mfu'_J\coloneqq\mfu'\cap J\mfu'$. By \eqref{eq:775}, we have $(\mfu')^{\perp_{\omega}}=\mfup\oplus \spa{X}$ for some $X\in \mfg\setminus \mfu$.

\begin{claim}
$JX\in \mfu'$.
\end{claim}

Pick $Y$ such that $\spa Y=\mfup$, so that $(\mfu')^{\perp_{\omega}}=\spa{X,Y}$. We need to show that $\omega(JX,X)=0=\omega(JX,Y)$. The first equality holds since $J$ is symmetric. For the second one, $\omega(JX,Y)=\omega(X,JY)=0$, since $X\in (\mfu')^{\perp_{\omega}}$ and $JY\in J\mfup\subset \mfu'$ by assumption.
\begin{flushright}
\qedsymbol{}
\end{flushright}

\medskip
Let us note that so
\[
\mfu'=J\mfup\oplus \mfu'_J \oplus \spa{JX}\,.
\]
Recall that $\mfu'\subset\mfg$ is a symplectic subspace, with symplectic form $\omega'$ induced by $\omega$. 
\begin{claim}\label{claim-omega}
$\omega'\big(\mfu'_J,J\mfup\oplus \spa{JX}\big)=0$.
\end{claim}

If $W\in J\mfup\oplus \spa{JX}$, there exist $Y\in \mfup$ and $\lambda\in \bR$ such that $W=JY+\lambda JX$. Thus, for any $Z\in \mfu'_J$, we have
\begin{equation*}
\omega'(Z,W)=\omega'(Z,JY)+\lambda \omega'(Z,JX)=\omega'(JZ,Y)+\lambda \omega'(JZ,X)=0
\end{equation*}
since $JZ\in \mfu'_J\subset\mfu'\subset \mfu$, while $Y\in \mfu^{\perp_{\omega}}$ and $X\in (\mfu')^{\perp_{\omega}}$ by assumption.
\begin{flushright}
\qedsymbol{}
\end{flushright}

\medskip
As a consequence, the restriction $\omega'_J$ of $\omega'$ to $\mfu'_J$ is non-degenerate as well. Now we have
\begin{align*}
\mfu &= \mfup\oplus\mfu'=\underbrace{\mfup\oplus J\mfup\oplus \mfu'_J }_{\mfu_J}\oplus\,\spa{JX}\,.
\end{align*}

This allows us to state the following result:
\begin{theorem}\label{th:complexsymplectic}
Let $\mfg$ be a $4n$-dimensional almost Abelian Lie algebra with a codimension one Abelian ideal $\mfu$, and let $(J,\omega)$ be an 
$\mathrm{Sp}(2n,\bC)$-structure on $\mfg$. Consider 
some complement $\mfu'$ of $\mfup$ in $\mfu$ which contains $J\mfup$ and set $\mfu'_J:=\mfu'\cap J\mfu'$. Take also $X\in \mfg\setminus \mfu$ such that $X\in (\mfu')^{\perp_{\omega}}$ and
\begin{equation*}
\mfu= \mfu'_J\oplus \mfup\oplus J\mfup\oplus \spa{JX}\,.
\end{equation*}
If we choose $Y\in \mfup\setminus \{0\}$ with $\omega(JY,JX)=1$, then
\begin{equation}\label{eq:splitting-vector-spaces}
(\mfg,J,\omega)=(\mfu'_J,J|_{\mfu'_J},\omega'_J)\oplus (V,J|_V,\omega_V)
\end{equation}
as complex symplectic vector spaces for
\begin{align*}
V&:=\mfup\oplus J\mfup\oplus \spa{JX,X}=\spa{Y,JY,JX,X},\\
\omega_V&:=J^*\alpha\wedge J^*\beta-\alpha\wedge \beta.
\end{align*}
and $\alpha,-J^*\alpha,\beta,-J^*\beta \in \mfg^*$ spanning the annihilator $\Ann{(\mfu_J')}\cong V^*$ of $\mfu'_J$ in $\mfg$ and being the dual basis of $(Y,JY,X,JX)$. Moreover, $(J,\omega)$ is a complex symplectic structure on $\mfg$ if and only if there exist $a,b,c\in \bR$, $u\in \mfu'_J$ and $f'_J\in \mathfrak{sp}(\mfu'_J,\omega'_J,J)$ such that
\begin{equation*}
f=\begin{pmatrix}
f'_J & \mathbf 0 & \mathbf 0 & u \\
\omega(Ju,\cdot) & a & 0 & b  \\
\omega(u,\cdot) & 0 & a & c  \\
\mathbf 0 & 0 & 0&  -a
\end{pmatrix}
\end{equation*}
with respect to the splitting $\mfu= \mfu'_J\oplus \mfup\oplus J\mfup\oplus \spa{JX} $ and to the basis $\{Y,JY,JX\}$ of $\mfup\oplus J\mfup\oplus \spa{JX}$.
\end{theorem}
\begin{proof}
Choose $0\neq Y\in\mfup$ such that $\omega(JY,JX)=1$. We first show that~\eqref{eq:splitting-vector-spaces} holds. Let us start observing that 
$\omega(Y,X)=-\omega(JY,JX)=-1$. As $JY\in \mfu'$ and $X\in (\mfu')^{\perp_{\omega}}$ by assumption, we also have $\omega(Y,JX)=\omega(JY,X)=0$. By Claim~\ref{claim-omega} one has $\omega'\big(\mfu'_J,J\mfu^{\perp_{\omega}}\oplus \spa{JX}\big)=0$,
and the symmetry of $J$, together with the $J$-invariance of $\mfu'_J$, yields
\begin{equation*}
\omega\big(\mfu'_J,\mfu^{\perp_{\omega}}\oplus \spa{X}\big)=0\,.
\end{equation*}
If we set
\begin{equation*}
V\coloneqq\mfup\oplus J\mfup\oplus \spa{JX,X}=\spa{Y,JY,JX,X}, 
\end{equation*}
and denote by $\alpha,\beta\in \mfg^*$ the elements in the annihilator of $\mfu'_J$ satisfying $\alpha(Y)=1=\beta(X)$, $\alpha(JY)=\alpha(JX)=\alpha(X)=0$ and $\beta(Y)=\beta(JY)=\beta(JX)=0$. We have
\begin{equation*}
(\mfg,J,\omega)=(\mfu'_J,J|_{\mfu'_J},\omega_J)\oplus (V,J|_V,\omega_V)
\end{equation*}
as complex symplectic vector spaces for
\begin{equation*}
\omega_V\coloneqq J^*\alpha\wedge J^*\beta-\alpha\wedge \beta.
\end{equation*}

For the second part of the statement, we set as usual $f\coloneqq\ad_X|_{\mfu}$. Then Theorem \ref{th:symplectic} implies the existence of some $a\in \bR$ with $f(Y)=a Y$ and Theorem \ref{th:complex} gives $f(JY)=Jf(Y)=a JY$. This shows how to define $f$ on $\mfup=\spa Y$ and $J\mfup=\spa{JY}$. Thus, it remains to describe $f$ on the elements of $\mfu'_J$ and on $JX$. For this, write $f|_{\mfu'}=f'+\gamma$ with $f'\in \mathfrak{gl}(\mfu')$ and $\gamma\in \mathrm{Hom}(\mfu',\mfup)\cong (\mfu')^*$. Note that by Theorem~\ref{th:symplectic}, we have $f'\in\mathfrak{sp}(\mfu',\omega')$. Also, by Theorem~\ref{th:complex} $f(\mfu_J')\subset f(\mfu_J)\subset \mfu_J$, so that $f'(\mfu'_J)\subset \mfu'_J\oplus \spa{JY}$. Consequently, there exist $\delta\in (\mfu'_J)^*$ and $f'_J\in \mathfrak{gl}(\mfu_J')$ with
\begin{equation*}
f'(Z)=\delta(Z) JY+ f'_J(Z)
\end{equation*}
for $Z\in \mfu'_J$. Besides, define $b,c\in \bR$ and $u\in \mfu'_J$ by
\begin{equation*}
f'(JX)=b JY+ c JX+ u.
\end{equation*}
Then $f'\in \mathfrak{sp}(\mfu',\omega')$ gives
\begin{equation*}
0=(f'.\omega')(JY,JX)=-\omega'(f'(JY),JX)-\omega'(JY,f'(JX))=-a-c,
\end{equation*}
i.e. $c=-a$. Moreover, for $Z\in \mfu_J'$, we get
\begin{equation*}
0=(f'.\omega')(Z,JX)=-\delta(Z)-\omega(Z,u),
\end{equation*}
which is equivalent to $\delta=\omega(u,\cdot)$. Hence, for any $Z\in \mfu_J'$, we have
\begin{equation*}
f(Z)=f'(Z)+\gamma(Z)Y=f_J'(Z)+\delta(Z)JY+\gamma(Z)Y
\end{equation*}
and $f(JZ)=Jf(Z)$ by Theorem \ref{th:complex}, which implies $\gamma(Z)=\delta(JZ)=\omega(u,JZ)=\omega(Ju,Z)$. Moreover, Theorem \ref{th:complex} and Theorem \ref{th:symplectic} together yield $f'_J\in \mathfrak{sp}(\mfu'_J,\omega'_J,J)\cong \mathfrak{sp}(2n-2,\bC)$. This gives the desired result.
\end{proof}

\begin{remark}
Note that the first part of the previous theorem is related to the complex symplectic reduction of a complex symplectic vector space by a coisotropic, $J$-invariant subspace. Indeed, we just proved that 
\[
\mfu'_J=\mfu_J/(\mfup+J\mfup)
\]
inherits a complex symplectic structure. $V$ is any complement of $\mfu'_J$ in $\mfg$, and it also inherits a complex symplectic structure. 
\end{remark}

\begin{remark}
In dimension 4, Theorem \ref{th:complexsymplectic} says that $\bR^3\rtimes_f \bR$ admits a complex symplectic structure if and only if there are $a,b,c\in \bR$ such that
\begin{equation*}
f=\begin{pmatrix}
a & 0 & b \\
0 & a & c  \\
0 & 0 & -a 
\end{pmatrix}
\end{equation*}
with respect to some basis of $\bR^3$. This agrees with the results of \cite[Theorem 3.1]{BGGL}, saying that the four-dimensional almost Abelian Lie algebras admitting a complex symplectic structure are $\mfr_{4,-1,-1}$, $\mfr \mfh_3$ and $\bR^4$. More precisely, we have:
\begin{itemize}
	\item if $a\neq 0$, then $\bR^3\rtimes_f \bR$ is isomorphic to $\mfr_{4,-1,-1}$;
	\item if $a=0$ but $(b,c)\neq (0,0)$, then $\bR^3\rtimes_f \bR$ is isomorphic to $\mfr \mfh_3$;
	\item if $a=b=c=0$, then $\bR^3\rtimes_f \bR$ is isomorphic to $\bR^4$.
\end{itemize}
\end{remark}
It is interesting to note that by \cite{BGGL}, a four-dimensional real almost Abelian Lie algebra has at most one complex symplectic structure, up to equivalence. We show that this is not true in general.
\begin{example}\label{ex:nonuniqueness}
We consider $(\mfg_i,J_0,\omega_0)$ with
\begin{equation*}
\begin{split}
J_0&=-e^1\otimes e_2+e^2\otimes e_1-e^3\otimes e_4+e^4\otimes e_3+e^5\otimes e_6-e^6\otimes e_5-e^7\otimes e_8+e^8\otimes e_7\\
\omega_0&=e^{14}+e^{23}-e^{58}+e^{67}
\end{split}
\end{equation*}
and $\mfg_i:=\bR^7\rtimes_{f_i} \bR$, $i=1,2$, where $f_1$ and $f_2$ are as in Theorem \ref{th:complexsymplectic} with $a_1=a_2=0$, $c_1=c_2=0$ and
\begin{equation*}
\begin{split}
u_1&=0\,,\\
u_2&=e_4\,,
\end{split}
\qquad
\begin{split}
b_1&=1\,,\\
b_2&=0\,,
\end{split}
\qquad 
\begin{split}
(f'_J)_1&=\left(\begin{smallmatrix} 0 & 0 & 1 & 0  \\ 0 & 0 & 0 & 1 \\
0 & 0 & 0 & 0\\ 0 & 0 & 0 & 0 \end{smallmatrix}\right)\,,\\
(f'_J)_2&=\mathbf 0\,.
\end{split}
\end{equation*}
By Theorem \ref{th:complexsymplectic}, $(\mfg_i,J_0,\omega_0)$ is an almost Abelian complex symplectic Lie algebra for $i=1,2$. Note that if $[\cdot,\cdot]_i$ denotes the Lie bracket of $\mfg_i$, the only non-zero Lie brackets of the standard basis $\{e_1,\ldots,e_8\}$ are given by
\begin{equation*}
\begin{split}
[e_8,e_3]_1&=e_1\,, \\
[e_8,e_1]_2&=-e_6\,, 
\end{split}
\qquad
\begin{split}
[e_8,e_4]_1&=e_2\,, \\
[e_8,e_2]_2&=-e_5\,, 
\end{split}
\qquad
\begin{split}
[e_8,e_7]_1&=e_5\,,\\
[e_8,e_7]_2&=e_4\,.
\end{split}
\end{equation*}
The Jordan normal forms of $f_1$ and $f_2$ are the same, hence $\mfg_1$ and $\mfg_2$ are isomorphic as Lie algebras. We claim that $(\mfg_1,J_0,\omega_0)$ and $(\mfg_2,J_0,\omega_0)$ are not isomorphic as complex symplectic Lie algebras. Assume on the contrary that such an isomorphism $F\colon(\mfg_1,J_0,\omega_0)\rightarrow (\mfg_2,J_0,\omega_0)$ exists. By Proposition \ref{pro:uniquenesscodim1ideal}, $\mfu=\spa{e_1,\ldots,e_7}$ is the unique Abelian ideal in $\frg_i$ for $i=1,2$, and $F$ has to preserve it. Hence, $F$ also preserves $\mfu^{\perp_{\omega_0}}=\spa{e_5}$.
Thus,
\begin{equation*}
F\big([e_8,e_7]_1\big)=F(e_5)=\lambda e_5
\end{equation*}
for some $\lambda\in \bR\setminus \{0\}$. Next, write
\begin{equation*}
F(e_7)=a e_7+v
\end{equation*}
for $a\in \bR$ and $v\in \spa{e_1,\ldots,e_6}$. Then
\begin{equation*}
F(e_8)=-F(J_0 e_7)=-J_0 F(e_7)=-J_0 a e_7-J_0 v=a e_8-J_0 v
\end{equation*}
since $F$ commutes with $J_0$. As $F(\mfu)\subseteq \mfu$ and $\mfu\oplus \spa{e_8}=\mfg_i$, $i=1,2$, and $F$ has to be surjective, we must have $a\neq 0$. But so
\begin{equation*}
\begin{split}
\lambda e_5& =F\big([e_8,e_7]_1\big)=[F(e_8),F(e_7)]_2=[a e_8-J_0 v, ae_7+v]_2=a^2 [e_8,e_7]_2+a [e_8,v]_2\\
&= a^2 e_4+ a [e_8,v]_2,
\end{split}
\end{equation*}
which is a contradiction since $a\neq 0$ and $[e_8,v]_2\in \spa{e_5,e_6}$ for $v\in \spa{e_1,\ldots,e_6}$. Thus, $(\mfg_1,J_0,\omega_0)$ and $(\mfg_2,J_0,\omega_0)$ are not isomorphic as complex symplectic Lie algebras and $\mfg\coloneqq\mfg_1$ admits two non-equivalent complex symplectic structures.

Note that $\mfg$ is a nilpotent Lie algebra with rational structure constants with respect to the chosen basis. Hence, the associated simply connected Lie group $G$ admits a lattice and so the compact almost Abelian nilmanifold $M\coloneqq\Gamma\backslash G$ admits two invariant complex symplectic structures which are not equivalent by diffeomorphisms of $M$ induced by Lie group automorphisms of $G$.
\end{example}

\section{Classification of complex symplectic almost Abelian Lie algebras}\label{sec:classification}

In this section we classify the almost Abelian Lie algebras $\mfg=\bR^{4n-1}\rtimes_f \bR$ admitting a complex symplectic structure $(J,\omega)$ in terms of the Jordan normal forms of $f$. We then
use this classification to provide conditions on the uniqueness of complex symplectic structures, up to equivalence, on almost Abelian Lie algebras - see Corollary \ref{co:uniqueness}. Moreover, we give explicit examples of compact almost Abelian solvmanifolds admitting complex symplectic structures in every possible dimension.

First, let us note that in Theorem \ref{th:complexsymplectic} one can always find a basis for the
Abelian subspace $\mfu'_J$ where the complex symplectic structure $(J|_{\mfu'_J},\omega'_J)$ takes the canonical form.
As a consequence, applying again the theorem if necessary, one may assume that for any complex symplectic almost Abelian Lie algebra $(\mfg,J,\omega)$ there exists a basis $\{e_k\}_{k=1}^{4n}$
for $\mfg$ where $(\mfg,J,\omega)=(\mathbb R^{4n-1}\rtimes_f \mathbb R,J_0,\omega_0)$ for
\begin{equation}\label{eq:formoff}
	f=\begin{pmatrix}
		A & \mathbf 0 & \mathbf 0 & u \\
		\omega_0(J_0 u,\cdot) & a & 0 & b  \\
		\omega_0(u,\cdot) & 0 & a & c  \\
		\mathbf 0 & 0 & 0&  -a
	\end{pmatrix}
\end{equation}
with $A\in \mathrm{sp}(2n-2,\mathbb C)\subset \mathbb R^{(4n-4)\times (4n-4)}$, $a,b,c\in \mathbb R$ and $u\in \mathbb R^{4n-4}$ together with
\begin{equation}\label{eq:J0omega0}
	\begin{split}
		J_0&=\sum_{k=1}^{2(n-1)} \left(-e^{2k-1}\otimes e_{2k}+e^{2k}\otimes e_{2k-1}\right)+e^{4n-3}\otimes e_{4n-2}-e^{4n-2}\otimes e_{4n-3}\\
		&-e^{4n-1}\otimes e_{4n}+e^{4n}\otimes e_{4n-1},\\
		\omega_0&=\sum_{l=1}^{n-1} \left(e^{4l-3}\wedge e^{4l}+e^{4l-2}\wedge e^{4l-1}\right)-e^{4n-3}\wedge e^{4n}+e^{4n-2}\wedge e^{4n-1}.
	\end{split}
\end{equation}
Note that the elements $e_{4n-3}$, $e_{4n-2}$, $e_{4n-1}$, $e_{4n}$ on the basis play the role of $(Y,J_0Y,J_0X,X)$. Indeed, for simplicity we will make use of the latter notation for them. 
From now on, we will assume that the complex symplectic structure $(J,\omega)$ on $\frg=\mathbb R^{4n-1}\rtimes_f\mathbb R$ takes the canonical form~\eqref{eq:J0omega0}. We want to provide the list of those $f$ that give rise to non-equivalent complex symplectic almost Abelian Lie algebras $(\mathbb R^{4n-1}\rtimes_f \mathbb R,J_0,\omega_0)$. To do so, we study the equivalence classes of complex symplectic almost Abelian Lie algebras.

\begin{lemma}\label{lemma-aux}
Assume that the almost Abelian Lie algebra $\frg=\mathbb R^{4n-1}\rtimes_f \mathbb R$ is not Abelian and also not isomorphic to $\mfh_3\oplus\bR^{4n-3}$ as a Lie algebra.
 Then any complex symplectic almost Abelian Lie algebra $(\tilde\mfg, \tilde J, \tilde \omega)=(\mathbb R^{4n-1}\rtimes_{\tilde f} \mathbb R,J_0,\omega_0)$ isomorphic to 
$(\mfg, J, \omega)=(\mathbb R^{4n-1}\rtimes_f \mathbb R,J_0,\omega_0)$
is given by
\[
\tilde f=\begin{pmatrix}
	\tilde A & \mathbf 0 & \mathbf 0 & \tilde u \\
	\omega_0(J_0\tilde u,\cdot) & \frac{a}{\lambda} & 0 & \tilde b  \\
	\omega_0(\tilde u,\cdot) & 0 & \frac{a}{\lambda} & \tilde c  \\
	\mathbf 0 & 0 & 0&  -\frac{a}{\lambda}
\end{pmatrix}
\]
with
\begin{equation*}
\begin{split}
\tilde A&=\frac{1}{\lambda}\Delta\,A\,\Delta^{-1},\\
\tilde u&=\frac{1}{\lambda^2}\big(\Delta u
	-a\,J\tilde u_X -\lambda\tilde A\,(J\tilde u_X) \big),\\
\tilde b&=\frac{1}{\lambda^3}\big(b+2\,\lambda\,a\,\mu_2
		-\omega(\tilde u_X,\,\Delta u+\lambda^2\tilde u)\big),\\
\tilde c&=\frac{1}{\lambda^3}\big( c-2\,\lambda\,a\,\mu_1
		+\omega( J\tilde u_X,\, \Delta u +\lambda^2\tilde u) \big)\,,
\end{split}
\end{equation*}
where $\lambda\in\mathbb R^*$, $\mu_1,\,\mu_2\in\mathbb R$, $\tilde u_X\in\tilde\mfu'_J$, and $\Delta\colon\mfu'_{J_0}\to\tilde\mfu'_{J_0}$ is an isomorphism 
of $4(n-1)$-dimensional Abelian Lie algebras such that 
\[
\Delta\,J_0\vert_{\mfu'_{J_0}}-J_0\vert_{\tilde\mfu'_{J_0}}\,\Delta=0
\quad \text{and}\quad
\Delta^*\,\omega_0\vert_{\tilde \mfu'_{J_0}}=\omega_0\vert_{\mfu'_{J_0}}\,.
\]
Moreover, any isomorphism between $(\mathbb R^{4n-1}\rtimes_f \mathbb R,J_0,\omega_0)$
and $(\mathbb R^{4n-1}\rtimes_{\tilde f} \mathbb R,J_0,\omega_0)$ preserving $(J_0,\omega_0)$
has the following form
\renewcommand{\arraystretch}{1.2}
\begin{equation*}
\tilde K=\left(\begin{array}{c|ccc|c}
\Delta & \mathbf 0 & \mathbf 0 & J\tilde u_{X} & \tilde u_X \\ \hline
-\frac{1}{\lambda}\,\omega(\tilde u_X,\Delta\,\cdot\,) & \frac{1}{\lambda} & 0 & -\mu_2 & \mu_1\\
\frac{1}{\lambda}\,\omega(J\tilde u_X,\Delta\,\cdot\,) & 0 & \frac{1}{\lambda} & \mu_1 & \mu_2\\
\mathbf 0 & 0 & 0 & \lambda & 0\\ \hline
\mathbf 0 & 0 & 0 & 0 & \lambda
\end{array}\right).
\end{equation*}
\end{lemma}

\begin{proof}
Let $(\mfg, J, \omega)=(\mathbb R^{4n-1}\rtimes_f \mathbb R,J_0,\omega_0)$ be a complex symplectic almost Abelian Lie algebra given by~\eqref{eq:formoff} and~\eqref{eq:J0omega0} in terms of a basis $\{e_k\}_{k=1}^{4n}$. We consider another complex symplectic almost Abelian Lie algebra $(\tilde\mfg, \tilde J, \tilde \omega)$ for which, by the same argument as above, there is a basis $\{v_k\}_{k=1}^{4n}$ for $\tilde\mfg$ where $(\tilde\mfg,\tilde J,\tilde \omega)=(\mathbb R^{4n-1}\rtimes_{\tilde f} \mathbb R,J_0,\omega_0)$ for
\begin{equation}\label{eq:formoff2}
	\tilde f=\begin{pmatrix}
		\tilde A & \mathbf 0 & \mathbf 0 & \tilde u \\
		\omega_0(J_0\tilde u,\cdot) & \tilde a & 0 & \tilde b  \\
		\omega_0(\tilde u,\cdot) & 0 & \tilde a & \tilde c  \\
		\mathbf 0 & 0 & 0&  -\tilde a
	\end{pmatrix}
\end{equation}
with $\tilde A\in \mathrm{sp}(2n-2,\mathbb C)\subset \mathbb R^{(4n-4)\times (4n-4)}$, 
$\tilde a,\tilde b,\tilde c\in \mathbb R$, $\tilde u\in \mathbb R^{4n-4}$, and 
$(J_0,\omega_0)$ as in~\eqref{eq:J0omega0} but written in terms of the $v$'s. Note that we will use $\{\tilde Y,J_0\tilde Y,J_0\tilde X,\tilde X\}$ instead of $\{v_{4n-3},v_{4n-2},v_{4n-1},v_{4n}\}$ and similarly for the corresponding elements of the basis $\{e_k\}_{k=1}^{4n}$.

\medskip
For simplicity, let us denote $(J,\omega)$ the expression of $(J_0,\omega_0)$ given in \eqref{eq:J0omega0}, both in terms of $e$'s and~$v$'s.

\medskip
Suppose $\mfg$, $\tilde\mfg$ are real almost Abelian Lie algebras with unique codimension one Abelian ideals
$\mfu$ and $\tilde\mfu$, respectively; by Proposition~\ref{pro:uniquenesscodim1ideal}, they are not isomorphic to $\mfh_3\oplus \bR^{4n-3}$. Then, any isomorphism $\phi\colon\mfg\to\tilde\mfg$ satisfies $\phi\,\mfu=\tilde\mfu$.
Consequently, any isomorphism from $\mfg=\mfu'_J\oplus \mfup\oplus J\mfup\oplus \spa{JX}\oplus\spa X$ to $\tilde\mfg=\tilde\mfu'_J\oplus \tilde\mfu^{\perp_\omega}\oplus J\tilde\mfu^{\perp_\omega}\oplus \spa{J\tilde X}\oplus\spa{\tilde X}$ is given, in terms of the bases $\{e_k\}_{k=1}^{4n}$ and $\{v_k\}_{k=1}^{4n}$ above, by
\begin{equation}\label{eq:auto}
\tilde K=\left(\begin{array}{c|ccc|c}
\Delta & \tilde u_Y & \tilde u_{JY} & \tilde u_{JX} & \tilde u_X \\ \hline
\delta_1 & r_1 & s_1 & t_1 & \mu_1\\
\delta_2 & r_2 & s_2 & t_2 & \mu_2\\
\delta_3 & r_3 & s_3 & t_3 & \mu_3\\ \hline
\mathbf 0 & 0 & 0 & 0 & \lambda
\end{array}\right),
\end{equation}
where $\lambda\in\mathbb R^*$, $r_i,\,s_i,\,t_i,\,\mu_i\in\mathbb R$ and 
$\delta_i:\mfu'_J\to\mathbb R$ for $i=1,2,3$, $\tilde u_X, \tilde u_Y, \tilde u_{JY}, \tilde u_{JX}\in\mathbb R^{4n-4}$,
and $\Delta$ the matrix associated to a linear map from $\mfu'_J$ to $\tilde\mfu'_J$. By abuse of notation, $\tilde u_X, \tilde u_Y, \tilde u_{JY}, \tilde u_{JX}$ will also denote the elements in $\tilde\mfu'_J$ and $\Delta\colon\mfu'_J\to \tilde\mfu'_J$. Let us observe that the identity $\phi\big([X,v]\big)=[\phi(X),v]$ for any $v\in \mfu$ yields the equality
\begin{equation}\label{eq:isomorphism}
\tilde K\,f-\lambda\,\tilde f\tilde K=0.
\end{equation}
for $\tilde K$ on $\mfu$. We are interested in those $\phi$ that preserve the complex symplectic structure $(\omega,J)$. Consequently, the following conditions prescribed by Definition~\ref{def:equivalence-CS} must be satisfied:
\begin{align}
\phi\circ J&=J\circ \phi\,,\label{eq:equivalence-J}\\
\phi^*\omega&=\omega\,.\label{eq:equivalence-omega}
\end{align}
Further, since $\phi$ is invertible, $\det \tilde K\neq 0$. Applying~\eqref{eq:equivalence-J} to an isomorphism $\phi$ given by~\eqref{eq:auto}
and to each element in $\mfg$, one
obtains
\begin{equation}\label{eq:conditions1}
\begin{split}
&\Delta J-J\Delta=0 \text{ \ on \,} \mfu'_J\\
&\delta_2=-\delta_1\circ J, \\
&\delta_3 = 0, \\ 
\quad
\end{split}
\qquad
\begin{split}
&\tilde u_{JY} = J\tilde u_Y,\\
&s_1 = -r_2, \\
&s_2 = r_1,\\
&s_3 = r_3 = \mu_3 =0,
\end{split}
\qquad\quad
\begin{split}
&\tilde u_{JX} = J\tilde u_X,\\
&t_1 = -\mu_2,\\
&t_2 = \mu_1,\\
&t_3 =\lambda.
\end{split}
\end{equation}
Now observe that for any $U_1,U_2\in \mfu'_J=\spa{e_1,\ldots,e_{4n-4}}$, \eqref{eq:equivalence-omega} gives
$\omega(U_1,\,U_2)=\omega(\Delta U_1,\,\Delta U_2)$, using~\eqref{eq:conditions1} and the fact that $(J,\omega)$ is a complex symplectic structure. As a consequence, 
$\Delta\colon\mfu'_J\to\tilde\mfu'_J$ is bijective. Indeed, suppose $\Delta$ 
is not injective. Then, there exists $0\neq U\in\mfu'_J$ such that $\Delta U=0$. 
As $\omega\vert_{\mfu'_J}$ is non-degenerate, there is $0\neq V\in\mfu'_J$ such that 
$0\neq\omega(U,V)=\omega(\Delta U,\,\Delta V)=0$, but this is a contradiction. 
Therefore, $\Delta$ is
injective and thus bijective, since $\text{dim}\,\mfu'_J=\text{dim}\,\tilde\mfu'_J$. 
Applying now~\eqref{eq:equivalence-omega} to the other possible pairs $(A,B)$ with 
$A,B\in\{e_k\}_{k=1}^{4n}$, one has
\begin{equation}\label{eq:conditions2}
\begin{split}
&\Delta \text{ bijective}, \\
&\delta_1= -\frac{1}{\lambda}\,\omega(\tilde u_X,\Delta\,\cdot\,),
\end{split}
\qquad
\begin{split}
&\Delta^*\omega=\omega \text{ \ on \,} \mfu'_J, \\
&r_1=\frac{1}{\lambda}, \\
\end{split}
\qquad
\begin{split}
&\tilde u_Y=0,\\[6pt]
&r_2=0.
\end{split}
\end{equation}

Finally, the condition \eqref{eq:isomorphism} for the reduced version of $\tilde K$ obtained
by considering~\eqref{eq:conditions1} 
and~\eqref{eq:conditions2} gives the desired result.
\end{proof}


Hence, to determine the possible real Jordan normal forms of $f$ for an almost Abelian Lie algebra $\bR^{4n-1}\rtimes_f \bR$, it is surely of importance to know the real Jordan normal forms
of $A\in \mathfrak{sp}(2n-2,\bC)$, considered as a real $(4n-4)\times (4n-4)$-matrix. For this purpose, we introduce the following notation:
\begin{notation}\label{not:blocks}
Let $F\in \End(\bR^N)$ and $m\in \bN$. We denote by $N_F(m,a)\in \bN_0$  the number of real Jordan blocks $J_m(a)$ of size $m$ for the eigenvalue $a\in \bR$ in the real Jordan normal form of $F$. If $z\in \bC\setminus \bR$, then $N_F(m,z)\in \bN_0$ will stand for the number of real Jordan blocks $J_m(z)$ of size $2m$ for the pair of complex conjugate eigenvalues $z=a+ib$, $\overline{z}=a-ib$, where we recall that for $M_z:=\left(\begin{smallmatrix} a & -b \\ b & a\end{smallmatrix}\right)$ we have
\renewcommand{\arraystretch}{1}
\begin{equation*}
J_m(z)=\begin{pmatrix}
         M_z & I_2 & & \\
         &  \ddots & \ddots & \\
         & & \ddots & I_2 \\
         & & & M_z 
        \end{pmatrix}\in \bR^{2m\times 2m}.
\end{equation*}
Similarly, for $G\in \End(\bC^N)$, $m\in \bN$ and $z\in \bC$, we denote by $N^{\bC}_G(m,z)$ the number of complex Jordan blocks $J_m(z)$ of size $m$ with $z$ on the diagonal.
\end{notation}

The possible real Jordan normal forms of $A$ follow directly from the possible complex Jordan normal forms of the complex $(2n-2)\times (2n-2)$-matrix $A$. Set:
\[D^+=\begin{pmatrix} 1 & 0 \\ 0 & 0 \end{pmatrix}, \qquad 
	D_-=\begin{pmatrix} 0 & 0 \\ 0 & -1 \end{pmatrix}, \qquad
	D(z)=\begin{pmatrix} z & 0 \\ 0 & -z \end{pmatrix}\,,
\]
where $z\in\mathbb C$. The following well-known result is of much use:
\begin{proposition}\label{pro:complexJNFsp2mC}
Let $\omega^{\bC}_0:=\sum_{i=1}^m e^{2i-1}\wedge e^{2i}$ denote the standard symplectic structure on the complex vector space $\bC^{2m}$ and let
\begin{equation*}
\mathfrak{sp}(2m,\bC):=\left\{\left. A\in \bC^{2m\times 2m}\right|\omega^{\bC}_0(Av,w)=-\omega^{\bC}_0(v,Aw)\textrm{ for all }v,w\in \bC^{2m}\right\}.
\end{equation*}
\begin{enumerate}
	\item[(a)]
	$A\in \bC^{2m\times 2m}$ is similar to a complex matrix in $\mathfrak{sp}(2m,\bC)$ if and only if for any $z\in \bC\setminus \{0\}$ and any $k\in \bN$ we have 
	$N^{\bC}_A(k,z)=N^{\bC}_A(k,-z)$ and 
	$N^{\bC}_A(2l-1,0)\equiv 0$ (mod $2$)
	for any $l\in \bN$.
	\item[(b)]
	Let $A\in\mathfrak{sp}(2m,\bC)$. Then $(\bC^{2m},\omega^{\bC}_0)=\sum_{i=1}^k (V_i,\omega_i)$ decomposes into a sum of symplectic $A$-invariant and $A$-irreducible subspaces of $\bC^{2m}$ with $\text{dim}\,V_i=2m_i$, and for each $(V_i,\omega_i)$, $i\in \{1,\ldots,k\}$, one of the following holds:
	\begin{itemize}
		\item[(i)]
		there exists some $z\in \bC\setminus \{0\}$ and a basis $v_1,\ldots,v_{2m_i}$ of $V_i$ such that $\omega_i=\sum_{j=1}^{m_i} v^{2j-1}\wedge v^{2j}$ and
		\begin{equation*}
		 A|_{V_i}=\begin{pmatrix}
		 D(z) &  D_- & & \\
		 D^+ & \ddots & \ddots & \\
		 & \ddots & \ddots &  D_- \\
		 & &  D^+ & D(z)
		 \end{pmatrix}
		 \end{equation*}
		 with respect to the ordered basis $(v_1,\ldots,v_{2m_i})$ of $V_i$,
		 \item[(ii)]
		 $m_i$ is odd, namely,
		 $m_i=2k-1$ for some $k\in \bN$, and there is a basis $v_1,\ldots,v_{4k-2}$ of $V_i$ 
		 such that 
		  $\omega_i=\sum_{j=1}^{2k-1} v^{2j-1}\wedge v^{2j}$ and
		  \begin{equation*}
		A|_{V_i}=\begin{pmatrix}
		0 &  D_- & & \\
		D^+ & \ddots & \ddots & \\
		& \ddots & \ddots &  D_- \\
		& &  D^+ & 0
		\end{pmatrix}
		\end{equation*}
		  with respect to the ordered basis $(v_1,\ldots,v_{4k-2})$ of $V_i$,
		   \item[(iii)]
		  there is a basis $v_1,\ldots,v_{2m_i}$ of $V_i$ such that 
		  $\omega_i=\sum_{j=1}^{m_i} v^{2j-1}\wedge v^{2j}$ and
		  \begin{equation*}
		  A|_{V_i}=\begin{pmatrix}
		  0 &  D_- & & \\
		  D^+ & \ddots & \ddots & \\
		  & \ddots & 0 &  D_- \\
		  & &  D^+ & N 
		  \end{pmatrix},
		  \end{equation*}
		  where $N:=\left(\begin{smallmatrix} 0 & 0\\ 1 & 0 \end{smallmatrix}\right)$,
		  with respect to the ordered basis $(v_1,\ldots,v_{2m_i})$ of~$V_i$.
	\end{itemize}
	\item[(c)]
	$A,B\in \mathfrak{sp}(2m,\bC)$ are similar to each other if and only if they are \emph{symplectically similar}, i.e. if there exists some $T\in \mathrm{Sp}(2m,\bC)$ with $B=TAT^{-1}$.
\end{enumerate}
\end{proposition}
\begin{proof}
Part (a) can be found for instance in \cite[Theorem 2.7]{MMRR} - see also the references mentioned there. Part (b) can be easily deduced from the explicit description of the normal forms of the matrix representing $\omega^{\bC}|_{V_i}$ in \cite[Theorem 2.7]{MMRR} by bringing these normal forms into the standard form on $V_i$.

For (c), note that $A,B\in \mathfrak{sp}(2m,\bC)$ are similar if and only if they have the same complex Jordan normal form and so the decompositions into $A$-invariant and $A$-irreducible symplectic subspaces in (b) is, up to a permutation, the same for all $A$ and $B$. Now observe that the map sending the standard basis to the concatenation of bases as in part (b) (i) $-$ (iii) is an element of $\mathrm{Sp}(2m,\bC)$ and so $A$ and $B$ are symplectically similar.
\end{proof}

Before we prove the classification and uniqueness results, we first show that one may further simplify the form of $f$ as in \eqref{eq:formoff}. For this, we set
\begin{equation*}
\begin{split}
I_{2,0}&:=\diag(1,1,0,0),\qquad I_{0,2}:=\diag(0,0,1,1),\\
D&:=\diag(-1,-1,1,1),\qquad \tilde{N}=\begin{pmatrix} 0 & 0 \\ I_2 & 0 \end{pmatrix}\in \bR^{4\times 4},\\
\end{split}
\end{equation*}

\begin{equation*}
\begin{split}
\tilde{J}_m(-1)&:=\begin{pmatrix}
D & -I_{0,2}  &  &  0\\
I_{2,0} & \ddots & \ddots &  \\
&  \ddots & D & -I_{0,2}\\
0 &  & I_{2,0} & D
\end{pmatrix}\in \bR^{4m\times 4m}, \\
\tilde{J}_{2k-1}&:=\begin{pmatrix}
0 & -I_{0,2}  &  &  0\\
I_{2,0} & \ddots & \ddots &  \\
&  \ddots & 0 & -I_{0,2}\\
0 &  & I_{2,0} & 0
\end{pmatrix}\in \bR^{(8k-4)\times (8k-4)},\\
\tilde{J}_{2k}&:=\begin{pmatrix}
0 & -I_{0,2}  &  &  0\\
I_{2,0} & \ddots & \ddots &  \\
&  \ddots & 0 & -I_{0,2}\\
0 &  & I_{2,0} & \tilde{N}
\end{pmatrix}\in \bR^{4k\times 4k}
\end{split}
\end{equation*}

\begin{proposition}\label{pro:classsimplification}
Let $(\mfg,J,\omega)$ be a $4n$-dimensional almost Abelian complex symplectic Lie algebra. Then $(\mfg,J,\omega)$ is isomorphic as a complex symplectic Lie algebra to $(\bR^{4n-1}\rtimes_f \bR,J_0,\omega_0)$ with $f$ being one of the following matrices
\begin{equation*}
\begin{split}
	&\begin{pmatrix}
	A & \mathbf 0 & \mathbf 0 & \mathbf 0 \\
	\mathbf 0 & 1 & 0 & 0  \\
	\mathbf 0 & 0 & 1 & 0  \\
	\mathbf 0 & 0 & 0&  -1
	\end{pmatrix},\  
	\begin{pmatrix}
	A & \mathbf 0 & \mathbf 0 & \mathbf 0 \\
	\mathbf 0 & 0 & 0 & b  \\
	\mathbf 0 & 0 & 0 & c  \\
	\mathbf 0 & 0 & 0&  0
	\end{pmatrix},\
		\begin{pmatrix}
	B & \mathbf 0 & \mathbf 0 & \mathbf 0 & \mathbf 0 \\
	\mathbf 0 & \tilde{J}_p(-1) & \mathbf 0 & \mathbf 0 & \mathbf{e}_{1}\\
	\mathbf 0 & - \mathbf{e}_{3}^T & 1 & 0 & 0\\
	\mathbf 0 &   \mathbf{e}_{4}^T & 0 & 1 & 0  \\
	\mathbf 0 & \mathbf 0 & 0&  0 & -1
	\end{pmatrix},\\
	&
	\begin{pmatrix}
	C & \mathbf 0 & \mathbf 0 & \mathbf 0 & \mathbf 0 \\
	\mathbf 0 & \tilde{J}_{2r-1} & \mathbf 0 & \mathbf 0 &  \mathbf{e}_{1}\\
	\mathbf 0 & - \mathbf{e}_{3}^T & 0 & 0 & b\\
	\mathbf 0 & \mathbf{e}_{4}^T  & 0 & 0 & c  \\
	\mathbf 0 & \mathbf 0 & 0&  0 & 0
	\end{pmatrix}, \ \begin{pmatrix}
	D & \mathbf 0 & \mathbf 0 & \mathbf 0 & \mathbf 0 \\
	\mathbf 0 & \tilde{J}_{2s} & \mathbf 0 & \mathbf 0 &  \mathbf{e}_1\\
	\mathbf 0 & - \mathbf{e}_3^T & 0 & 0 & b\\
	\mathbf 0 &  \mathbf{e}_4^T & 0 & 0 & c \\
	\mathbf 0 & \mathbf 0 & 0&  0 & 0
	\end{pmatrix},
	\end{split}
	\end{equation*}
	for some $b,c\in \bR$ and $A\in \mathfrak{sp}\big(2n-2,\bC\big)$,  
	$B\in \mathfrak{sp}\big(2(n-1-p),\bC\big)$ with $p\in \{1,\ldots,n-1\}$, 
	$C\in \mathfrak{sp}\big(2(n-2r),\bC\big)$ for some $r\in \{1,\ldots,\lfloor\tfrac{n}{2}\rfloor\}$, 
	$D\in \mathfrak{sp}\big(2(n-1-s),\bC\big)$ for $s\in \{1,\ldots,n-1\} $, seen as real matrices of double size.
	The elements $\mathbf{e}_k$ stand for the vectors of the canonical basis of 
	$\mathbb R^N$ with $N\in\{4p, 8r-4, 4s\}$, depending on the case.
\end{proposition}

\begin{proof}
First consider the case where $\mfg$ is isomorphic to $\mfh_3\oplus \bR^{4n-3}$ as a Lie algebra.
We show that then $(\mfg,J,\omega)\cong (\bR^{4n-4},J_1,\omega_1)\oplus (\mfh_3\oplus \bR,J_2,\omega_2)$
as complex symplectic Lie algebras, where $(J_1,\omega_1)$ is a complex symplectic structure on $\bR^{4n-4}$ and $(J_2,\omega_2)$ is a complex symplectic structure on $\mfh_3\oplus \bR$. 

To prove the aforementioned statement, take some $X\in \mfg\setminus \mfu$, where $\mfu$ is some codimension one Abelian ideal with $JX\in\mathfrak u$. As $\mfg=\mfh_3\oplus \bR^{4n-3}$, we know that $f:=\ad(X)|_{\mfu}$ satisfies $\mathrm{im}(f)\subseteq\ker(f)$ and $\mathrm{im}(f)$ is one-dimensional whereas $\ker(f)$ is $(4n-2)$-dimensional. Moreover, the vanishing of the Nijenhuis tensor of $J$ yields that $\ker(f)$ is $J$-invariant and so $\ker(f)=\mfu_J$. Since $\mfu=\mfu_J\oplus \spa{JX}$, we thus have $Y:=f(JX)\neq 0$. Set now $V:=\spa{X,JX,Y,JY}$. Then $V$ is a $J$-invariant ideal in $\mfg$ isomorphic to $\mfh_3\oplus \bR$. Moreover,
\begin{equation*}
0=d\omega(X,JX,U)=-\omega(f(JX),U)-\omega(JX,f(U))=-\omega(Y,U)
\end{equation*}
for any $U\in \ker(f)$. Consequently, the non-degeneracy of $\omega$ and $\mfg=\ker(f)\oplus \spa{X,JX}$ implies that $-\omega(JY,JX)=\omega(Y,X)\neq 0$ or $\omega(JY,X)=\omega(Y,JX)\neq 0$. Thus, $V$ is also $\omega$-symplectic. 
We then have
\begin{equation*}
(\mfg,J,\omega)\cong (V^{\perp \omega},\,J|_{V^{\perp \omega}},\,\omega|_{V^{\perp \omega}})\oplus (V,\,J|_V,\,\omega|_V)
\end{equation*}
as complex symplectic vector spaces. To prove the claimed result, we need to show that $V^{\perp \omega}$ is an Abelian ideal,
but this will come as a consequence of $V^{\perp \omega}\subset \ker(f)$. 
To prove the latter assumption,
let $U\in V^{\perp \omega}$ and write $U=\tilde{U}+\hat{U}$ for $\tilde{U}\in \ker(f)$ and $\hat{U}\in \spa{X,JX}$. 
As $Y, JY\in V$
and we have shown above that $\omega\big(\ker(f),Y\big)=\{0\}$, we get
\begin{equation*}
\begin{split}
0&=\omega(U,Y)=\omega(\hat{U},Y),\\ 0&=\omega(U,JY)=\omega(\tilde{U},JY)+\omega(\hat{U},JY)
	=\omega(J\tilde{U},Y)+\omega(\hat{U},JY)=\omega(\hat{U},JY).
\end{split}
\end{equation*}
Moreover, $\omega(\hat{U},X)=\omega(\hat{U},JX)=0$, and so the non-degeneracy of $\omega$ on $V=\spa{X,JX,Y,JY}$ yields $\hat{U}=0$, i.e. $U=\tilde{U}\in \ker(f)$. 
This gives the desired result.

Now, note that there is obviously only one complex symplectic structure on $\bR^{4n-4}$ up to equivalence, and the same is true for $\mfh_3\oplus \bR$ by \cite[Proposition 5.4]{BFLM}. Thus, also $\mfh_3\oplus \bR^{4n-3}$ admits only one complex symplectic structure up to equivalence and it can be obtained by the second matrix in the statement with $A=0$, $b=1$ and $c=0$.


Thus, we may now restrict to the case that $\mfg$ has a unique codimension one Abelian ideal $\mfu$. Note that by Theorem \ref{th:complexsymplectic}, we may assume that $(\mfg,J,\omega)=(\bR^{4n-1}\rtimes_f \bR, J_0,\omega_0)$ with $f$ as in \eqref{eq:formoff}. Now we try to simplify $f$ while keeping $(J_0,\omega_0)$ fixed. For this, we
will make use of Lemma~\ref{lemma-aux}, which gives us the type of changes $\tilde K$ that we
can apply to $(\bR^{4n-1}\rtimes_f \bR, J_0,\omega_0)$ to get an equivalent complex symplectic
almost Abelian Lie algebra $(\bR^{4n-1}\rtimes_{\tilde f} \bR, J_0,\omega_0)$ with possible simpler $\tilde f$.

We first observe that one can assume that either $u=\mathbf 0$ or $u$ lies in the generalised eigenspace of $A$ with eigenvalue $-a$ but not in the image of $A+a\,\mathbf I_{4n-4}$. To check this assumption, let us consider an isomorphism $\tilde K$ as in Lemma~\ref{lemma-aux} with
$$\Delta=\mathbf I_{4n-4}, \quad \tilde u_X=-Jv, \quad 
\lambda=1, \quad \mu_1=\mu_2=0,$$
where $v\in\mathbb R^{4(n-1)}$ is to be determined. Then, note that the element $\tilde u$ in the matrix $\tilde f$ is given by
$$\tilde u=u-(A+a\,\mathbf I_{4n-4})v.$$

If $-a$ is not an eigenvalue of $A$ then $\text{det}(A+a\,\mathbf I_{4n-4})\neq 0$ and one can choose $v=(A+a\,\mathbf I_{4n-4})^{-1}\,u$ to get $\tilde u=\mathbf 0$. 

If $-a$ is an eigenvalue of $A$ one can proceed as follows. Let us note that the space
$\mathbb R^{4(n-1)}$ 
can be decomposed into the sum of (maximal) generalised eigenspaces 
$\mathbb R^{4(n-1)} =\oplus_{\lambda \in \text{Spec}(A)} V_{\lambda}(A)$. Consequently, 
the vector $u\in\mathbb R^{4(n-1)}$ can be written as
$u=\oplus_{\lambda \in \text{Spec}(A)} u_{\lambda}$, where $u_{\lambda}\in V_{\lambda}(A)$ 
for each $\lambda \in \text{Spec}(A)$. Recall that each $V_{\lambda}(A)$ is invariant under $A$ 
and there is a basis of $V_{\lambda}(A)$ such that $A|_{V_{\lambda}(A)}$ is block-diagonal with Jordan blocks having $\lambda$ on the diagonal. Therefore, one has that 
$(A+a\,\mathbf I_{4n-4})|_{V_{\lambda}(A)}$ maps $V_{\lambda}(A)$ again into $V_{\lambda}(A)$ 
and, in terms of the previous basis, $(A+a\,\mathbf I_{4n-4})|_{V_{\lambda}(A)}$ is also block-diagonal but with Jordan blocks having $\lambda+a$ on the diagonal. For every $\lambda\neq -a$, one can then find $v_{\lambda}\in V_{\lambda}(A)$ such that 
$(A+a\,\mathbf I_{4n-4})v_{\lambda}=-u_{\lambda}$. Thus, setting 
$v=\oplus_{\lambda\in\text{Spec}(A)\setminus\{-a\}} v_{\lambda}$, we have that 
$\tilde{u}:=(A+a\,\mathbf I_{4n-4})v+u$ is in $V_{-a}(A)$. Moreover, we can assume that 
$\tilde u$ is not in the image of $A+a\,\mathbf I_{4n-4}$. Indeed, if we suppose that 
$\tilde{u}$ is in the image of $A+a\,\mathbf I_{4n-4}$, then there should exist some 
$w$ such that $(A+a\,\mathbf I_{4n-4})\,w=\tilde{u}$; applying a new change of basis as 
above one gets 
$\tilde{\tilde u}=\tilde u - (A+a\,\mathbf I_{4n-4})\,\tilde v=(A+a\,\mathbf I_{4n-4})(w-\tilde v)$ 
so taking $\tilde v=w$ gives $\tilde{\tilde u}=\mathbf 0$.

Now, if $a\neq 0$ in~\eqref{eq:formoff}, observe that Lemma~\ref{lemma-aux} can be applied
with $\tilde K$ defined by the following choices of the parameters:
$$\Delta=\mathbf I_{4n-4}, \quad \tilde u_X=\mathbf 0, \quad 
\lambda=a, \quad \mu_1=\frac{c}{2a^2}, \quad \mu_2=-\frac{b}{2a^2}.$$
Then, the entries of $\tilde f$ are given by
$$\tilde A=\frac{1}{a}\, A, \qquad \tilde u=\frac{1}{a}\,u, \qquad \tilde b=\tilde c=0.$$
In particular, note that $\frac{a}{\lambda}$ becomes equal to $1$ and if $u=\mathbf 0$ then also 
$\tilde u=\mathbf 0$. Moreover, let us remark that since $u$ lies in the generalised eigenspace of $A$ with eigenvalue $-a$, there exists some $k\in\mathbb N$ such that
$$0=(A+a\,\mathbf I_{4n-4})^k\,u=(a\,\tilde A+a\,\mathbf I_{4n-4})^k\,a\tilde u=
a^{k+1}(\tilde A+\mathbf I_{4n-4})^k\,\tilde u,$$
and thus $\tilde u$ belongs to the generalized eigenspace of $\tilde A$ with eigenvalue $-1$. In addition, if we had $\tilde u=(\tilde A+\mathbf I_{4n-4})\,\mathbf w$ for some 
$\mathbf w\in\mathbb R^{4n-4}$, then $u$ would be in the image of $A+a\,\mathbf I_{4n-4}$,
which is a contradiction.
As a consequence, one may 
assume without loss of generality that the initial matrix $f$ given by~\eqref{eq:formoff} has
\begin{equation}\label{eq:values-in-reduced-f-1}
(a,b,c)=\begin{cases}
(1,0,0), \text{ or}\\
(0,b, c), \text{ with }b,c\in\mathbb R,
\end{cases}
\end{equation}
and $u$ either equal to zero or in the generalised eigenspace of $A$ with eigenvalue $-a$ but not in the image of $A+a\,\mathbf I_{4n-4}$, for the previous two values of $a$. 
Moreover, observe that these two sets of tuples give rise to non-equivalent complex symplectic almost Abelian Lie algebras. We next distinguish two cases depending on the value of $u$.

If $u=\mathbf 0$, then one obtains the first two claimed forms for $f$ directly 
from the two choices~\eqref{eq:values-in-reduced-f-1} above, respectively.

If $u\neq\mathbf 0$, then this vector lies in the generalised eigenspace of $A$ with eigenvalue $-a$ and not in the image of $A+a\,\mathbf I_{4n-4}$, where $a\in\{0,1\}$ in view of~\eqref{eq:values-in-reduced-f-1}. Consequently, the complex Jordan chain of $u$ generates a complex Jordan block of $A$ of some size $m$ with $-a$ on the diagonal. Two possibilities arise.

Suppose that $a=1$. Then, one easily deduces from Proposition~\ref{pro:complexJNFsp2mC}~(a) 
that there is another complex Jordan block of the same size $m$ but with $1$ on the diagonal. 
Moreover, by Proposition~\ref{pro:complexJNFsp2mC}~(b) the complex generalised eigenvectors 
corresponding to these two Jordan blocks generate a symplectic $A$-invariant and $A$-irreducible 
space $V$ whose complement in $\mathbb C^{2n-2}$ is $A$-invariant symplectic. In fact, it is 
possible to find a basis $(w_1,\ldots,w_{2m})$ of $V$ where $\omega_0|_V$ is preserved 
and $A|_V$ is given by Proposition~\ref{pro:complexJNFsp2mC}~(b)~(i). Thus, considering 
now all complex matrices as real matrices of double size, we can take $\tilde K$ from 
Lemma~\ref{lemma-aux} with $\lambda=1$, $\tilde u_X=\mathbf 0$, $\mu_1=\mu_2=0$ and
$\Delta\in\mathrm{Sp}(2n-2,\bC)$ to be the matrix that brings $A$ into 
$$\tilde A=\Delta\,A\,\Delta^{-1}=\begin{pmatrix} B & 0 \\ 0 & \tilde{J}_m(-1) \end{pmatrix},$$
where $B\in \mathfrak{sp}\big(2(n-1-m),\bC\big)$, and also takes $u$ into the first element of $\{w_k\}_{k=1}^{2m}$ seen as a real basis (recall that $u$ is in the generalised eigenspace of $A$ with eigenvalue $-1$ but not in the image of $A-\mathbf I_{4n-4}$). In particular,
note that the real counterpart of the complex basis $\{w_k\}_{k=1}^{2m}$ precisely 
coincides with the elements $v_{4n-4m-3},\ldots,v_{4n-4}$ of the basis $\{v_k\}_{k=1}^{4n}$ in terms
of which the matrix $\tilde f$ is written. Hence, this allows to define $\tilde u$ and fix the 
remaining terms in $\tilde f$, namely,
\begin{equation*}
\begin{split}
\tilde u&=\Delta u=v_{4n-4m-3}\\
\omega_0(J_0\tilde u,\cdot) &=-\omega_0(v_{4(n-m)-2},\cdot)=-v^{4(n-m)-1}\\
\omega_0(\tilde u,\cdot) &=\omega(v_{4(n-m)-3},\cdot)=v^{4(n-m)}.
\end{split}
\end{equation*}
This gives the third claimed possible form of $f$ in the assertion.

Suppose now that $a=0$ and that the size $m$ of the Jordan block corresponding to $u$ is odd, namely,
$m=2l-1$. As a consequence of Proposition~\ref{pro:complexJNFsp2mC}~(a), there is
another complex Jordan block of the same size $2l-1$ with $0$ on the diagonal. Reasoning in
a similar way to previous case, one can find a basis $(w_k)_{k=1}^{2(2l-1)}$ for the 
space $V$ spanned by the generalised eigenvectors corresponding to 
these two Jordan blocks where $\omega_0|_V$ is preserved 
and $A|_V$ is given by Proposition~\ref{pro:complexJNFsp2mC}~(b)~(ii).
Making again all complex matrices to be real matrices of double size and following the same
ideas as above, it is possible to bring $f$ into the fourth claimed possible form in the assertion.

The remaining case is $a=0$ with an even size $m$ of the corresponding Jordan block, namely,
$m=2l$. A similar argument as in the two cases above, now following 
Proposition \ref{pro:complexJNFsp2mC}~(b)~(iii) for the aforementioned Jordan block, 
gives the fifth form for $f$ in the statement.
\end{proof}
Let us recall that a Lie algebra $\mfg$ is said to be \emph{unimodular} if $\text{tr}(\ad_X)=0$ for all $X\in\mfg$. In our complex symplectic almost Abelian setting, note that this is equivalent 
to $\text{tr}(f)=0$.
As a consequence of Proposition \ref{pro:classsimplification}, one gets the following result:
\begin{theorem}\label{th:almostAbelianLAclass}
An almost Abelian Lie algebra $\mfg=\bR^{4n-1}\rtimes_f \bR$ admits a complex symplectic structure if and only if 
\begin{equation*}
N_f(m,z)=N_f(m,-z),\quad N_f(m,ib)\equiv 0 \text{ (mod $2$)}
\end{equation*}
for all $z\in \bC\setminus (\bR\cup i\bR)$, all $b\in \bR$
and all $m\in \bN$, and one of the set of conditions in (a) (i), (a)  (ii) or (b) (i) -- (b) (iv) is satisfied. These conditions are as follows:
\begin{itemize}
	\item[(a)]
	if $\mfg$ is not unimodular, then for any $k,m\in\mathbb N$ one has
	\begin{equation*}
	\begin{split}
 & N_f(2k,0)\equiv 0 \text{ (mod $2$)},\\
 & N_f(2k-1,0)\equiv 0 \text{ (mod $4$)},\\
 & N_f(m,a_0)\equiv 0 \text{ (mod $2$)} \text{ for some }a_0\in \bR\setminus \{0\},\\
 & N_f(m,a)=N_f(m,-a),\ N_f(m,a)\equiv 0 \text{ (mod $2$)} \text{ for every } a\in \bR\backslash \{a_0,0,-a_0\};
\end{split}
\end{equation*}
moreover, one of the following holds:

\smallskip
\begin{itemize}
	\item[(i)] $N_f(1,a_0)=1+N_f(-1,a_0)$ and 
	$N_f(m,a_0)=N_f(m,-a_0)$, for any $m\geq 2$; or
	\item[(ii)]
	there exists some $m_0\in \bN$ such that $N_f(m_0,a_0)=N_f(m_0,-a_0)-1$, 
	$N_f(m_0+1,a_0)=N_f(m_0+1,-a_0)+1$, and $N_f(m,a_0)=N_f(m,-a_0)$, for
	any $m\in \bN\setminus \{m_0,m_0+1\}$.
\end{itemize}

\smallskip
\item[(b)]
if $\mfg$ is unimodular, then there exists some $k_0\in \bN$ such that for every $k\in \bN\setminus \{k_0\}$ and any $m\in\mathbb N$ one has
\begin{equation*}
\begin{split}
& N_f(2k,0)\equiv 0 \text{ (mod $2$)},\\
&N_f(2k-1,0)\equiv 0 \text{ (mod $4$)},\\
&N_f(m,a)=N_f(m,-a),\ N_f(m,a)\equiv 0 \text{ (mod $2$)}, \text{ for every }a\in \bR\setminus \{0\};
\end{split}
\end{equation*}
moreover, one of the following holds:

\smallskip
\begin{itemize}
	\item[(i)]
	$k_0=1$, \ $N_f(1,0)\equiv 3$ (mod $4$) \ and \ $N_f(2,0)\equiv 0$ (mod $2$);
	\item[(ii)]
	$k_0=1$, \ $N_f(1,0)\equiv 1$ (mod $4$) \ and \ $N_f(2,0)\equiv 1$ (mod $2$);

	\item[(iii)]
	$N_f(2k_0-1,0)\equiv 1$ (mod $4$) \ and \ $N_f(2k_0,0)\equiv 3$ (mod $4$);

	\item[(iv)] $k_0\geq 2$, \ 
	$N_f(2k_0-1,0)\equiv 1$ (mod $4$) \ and \ $N_f(2k_0,0)\equiv 1$ (mod $2$).
\end{itemize}
\end{itemize}
\end{theorem}
\begin{proof}
We first recall the following well-known facts from Linear Algebra. Let $M\in \bC^{N\times N}$ be a complex square matrix and let $v$ be a generalised eigenvector of $M$ with complex eigenvalue $z\in \bC$ generating a Jordan chain of length $m$. Consider now $M$ as a real square matrix $M\in \bR^{2N\times 2N}$ of double size and denote the multiplication by $i$ on $\bR^{2N}\cong \bC^N$ by $J$. Then, if $z=a\in \bR$, both $v$ and $iv$ are linearly independent generalised eigenvectors generating two different Jordan chains with eigenvalue $a$ of the same length $m$, and so $N_M(m,a)\equiv 0$ (mod 2). Moreover, if $z\in \bC\setminus \bR$, and we extend now $M$ to a complex-linear map on $\bC^{2N}$, then $v-iJv$ generates a complex Jordan chain with eigenvalue $z$ of length $m$ and $v+iJv$ generates a complex Jordan with eigenvalue $\overline{z}$ of length $m$. Together, these complex Jordan blocks give rise to a real Jordan block $J_m(z)$ of even dimension as shown in Notation \ref{not:blocks}.

Suppose now that $M\in\mathfrak{sp} (2l,\mathbb C)$. 
By Proposition~\ref{pro:complexJNFsp2mC}~(a) and the observation above, when $M$ is 
considered as a real $(4l\times 4l)$-matrix, one has
$$N_M(2k,0)\equiv 0 \ (mod \ 2) \text{ \ and \ }
N_M(2k-1,0)\equiv 0 \ (mod \ 4), \text{ \ for all \ } k\in\mathbb N.$$
Furthermore, for every $z\in\mathbb C\setminus\{0\}$, we get $N_M(m,z)=N_M(m,-z)$. In particular,
$$\begin{array}{ll}
\text{ for }z=a\in\mathbb R\setminus\{0\}:&  N_M(m,a)=N_M(m,-a), \ N_M(m,a)\equiv 0 \text{ (mod 2)},\\[2pt]
\text{ for }z=ib\in i\,\mathbb R\setminus\{0\}:&  N_M(m,ib)\equiv 0 \text{ (mod 2)},
\end{array}$$
where the last assertion follows from the fact that 
$\bar z=-ib=-z$ and this enables to combine $J_m(z)$ and $J_m(-z)$ appropriately.

Note that, given an almost Abelian Lie algebra $\mfg=\bR^{4n-1}\rtimes_f \bR$ with
complex symplectic structure, $f$ may be assumed to be, up to a non-zero scaling, as in 
Proposition~\ref{pro:classsimplification}. Hence, $N_f(m,z)=N_M(m,z)$ for every $m\in\mathbb N$ 
and every $z\in\mathbb C\setminus\mathbb R$, where $M=A,B,C$ or $D$ depending on the
matrix that represents~$f$. This fact together with the observations above concerning 
$M\in\mathfrak{sp} (2l,\mathbb C)$ give the first part of 
the statement.

We now need to separately study the real eigenvalues for each possible form of~$f$ in 
Proposition~\ref{pro:classsimplification}, as one has to combine the information coming from the
two diagonal blocks that conform $f$.

Let us first suppose that $\mfg$ is non unimodular. This corresponds to the first and third possible forms of $f$ in Proposition~\ref{pro:classsimplification}. Then, it is clear that for any $k\in\mathbb N$
one has $N_f(2k,0)=N_M(2k,0)\equiv 0$ (mod 2) and $N_f(2k-1,0)=N_M(2k-1,0)\equiv 0$ (mod 4),
where $M=A$ if the first form of $f$ holds and $M=B$ if the third one is considered. In addition,
$N_f(m,c)=N_M(m,c)=N_M(m,-c)=N_f(m,-c)\equiv 0$ (mod 2) for every $m\in\mathbb N$ and every
$c\in\mathbb R\setminus\{-1,\,0,\,1\}$. We separately study the eigenvalues $\pm1$ for each form
of $f$:

\smallskip\noindent
$\bullet$ \ Assume $f$ is of the first possible form in Proposition~\ref{pro:classsimplification}. Then,
$N_f(1,1)=N_A(1,1)+2$ and $N_f(1,-1)=N_A(1,-1)+1$. Concerning those blocks of size $m\geq 2$,
it is straightforward to see that $N_f(m,1)=N_A(m,1)=N_A(m,-1)=N_f(m,-1)\equiv 0$ (mod 2). This
gives case (a) (i).

\smallskip\noindent
$\bullet$ \ Let $f$ correspond to the third from in Proposition~\ref{pro:classsimplification}. Observe
that the matrix $\tilde J_p(-1)$ has two Jordan blocks of order $p$ for the eigenvalue $-1$ and 
two Jordan blocks of the same order for $1$. In particular, the eigenvectors associated to 
$1$ correspond
to the third and fourth elements of the basis with respect to which $\tilde J_p(-1)$ is written, and we will
denote them by $u, v$. Those for the eigenvalue $-1$ are the fourth and third last elements
of the same basis, and we will refer to them as $w, t$. To count the Jordan blocks of $f$
one needs to take into account how the entries $\mathbf e_1$, $-\mathbf e_3^T$ and 
$\mathbf e_4^T$
of the matrix $f$ interfere with the Jordan blocks of $\tilde J_p(-1)$. We first remark that $u,v$ are not
eigenvectors of $f$, but they lie in the generalized eigenspaces of the third and second last vectors
of the basis of $f$, respectively. This gives two Jordan blocks of size $p+1$ with $1$ on the diagonal. In particular,
$N_f(p+1,1)=N_B(p+1,1)+2$, $N_f(p,1)=N_B(p,1)-2$ and
 $N_f(m,1)=N_B(m,1)$ for every $m\in\mathbb N\setminus\{p, p+1\}$.
In contrast, $w,t$ are eigenvectors of $f$ but the order of the Jordan block associated to one of them 
increases its order with respect to $\tilde J_p(-1)$, 
due to the last column of $f$. Hence, $N_f(p+1,-1)=N_B(p+1,-1)+1$, \ 
$N_f(p,-1)=N_B(p,-1)-1$ and $N_f(m,-1)=N_B(m,-1)$ for $m\in\mathbb N\setminus\{p, p+1\}$.
It suffices to recall that $N_B(m,1)=N_B(m,-1)\equiv 0$ (mod 2) to get 
part (a) (ii) of the statement, simply renaming $p$.

\medskip
We next consider the unimodular case, so one may assume that $f$ is of the second, fourth or fifth form in Proposition \ref{pro:classsimplification}. We first observe that for any 
$a\in\mathbb R\setminus\{0\}$ one has
$N_f(m,a)=N_M(m,a)=N_M(m,-a)=N_f(m,-a)\equiv 0 \text{ (mod 2)},$
for every $m\in\mathbb N$, where $M=A, C$ or $D$ depending on the three possible forms of $f$.
We need to study in detail the Jordan blocks corresponding to the zero eigenvalue, as their order
and number will depend on the two diagonal blocks of $f$.

\smallskip\noindent
$\bullet$ \ Let $f$ be given by the second matrix in Proposition \ref{pro:classsimplification}. Note that
if $b=c=0$, then $N_f(1,0)=N_A(1,0)+3$ and $N_f(m,0)=N_A(m,0)$ for every 
$m\in\mathbb N\setminus\{1\}$. This gives part (b) (i) of the statement. Furthermore, if one has
$bc=0$ but $(b,c)\neq (0,0)$, then $N_f(1,0)=N_A(1,0)+1$ and $N_f(2,0)=N_A(2,0)+1$, from where
we get (b) (ii). Finally, observe that the case $bc\neq 0$ can be reduced to $\tilde b\tilde c=0$
applying a similarity transformation.

\smallskip\noindent
$\bullet$ \ Suppose $f$ has the fourth form in Proposition \ref{pro:classsimplification}. Then, 
one easily sees that one may bring $f$ by a similarity transformation into a form with 
$\tilde{b}=\tilde{c}=0$. Note that $\tilde J_{2r-1}$ gives four Jordan blocks of equal size $2r-1$.
To count the Jordan blocks of $f$ we need to study how the entries $\mathbf e_1$, $-\mathbf e_3^T$
and $\mathbf e_4^T$ interfere with the blocks of $\tilde J_{2r-1}$. Observe that the eigenvectors
of $\tilde J_{2r-1}$ correspond to the third and four elements of the basis in which this matrix
is written, now denoted by $u, v$, and the four and third last elements of the same basis,
named $w,t$. Observe that $u,v$ are not eigenvectors of $f$, but they lie in the generalized 
eigenspaces generated by the third and second last columns of $f$, respectively. Thus, $f$ has at
least two Jordan blocks associated to the eigenvalue $0$ of order $2r$ more than $B$. Furthermore,
$v$ and $w$ are eigenvectors of $f$, but the order of the Jordan block coming from $\tilde J_{2r-1}$
is only preserved for $t$ when seen in $f$. The one corresponding to $w$ 
has now order $2r$, due to the last column of the matrix $f$. Consequently, 
$N_f(2r,0)=N_C(2r,0)+3$ and
$N_f(2r-1,0)=N_C(2r-1,0)-3$. This gives part (b) (iii) of the theorem for $k_0:=r$.

\smallskip\noindent
$\bullet$ \ Assume now that $f$ is determined by the last matrix in Proposition \ref{pro:classsimplification}.
Without loss of generality, we can take $b=c=0$ as in the previous case. Let us remark that the 
matrix $\tilde J_{2s}$ has two Jordan blocks of size $2s$ and the eigenvectors that originate
these blocks precisely correspond to the third and fourth vectors of the basis where $\tilde J_{2s}$
is given. However, these are not eigenvalues of $f$ due to the entries $\mathbf e_1$, 
$-\mathbf e_3^T$ and $\mathbf e_4^T$ in the corresponding matrix. By a similar argument as above,
one can nonetheless check that $N_f(2s+1,0)=N_D(2s+1,0)+1$, $N_f(2s+2,0)=N_D(2s+2,0)+1$ and $N_f(2s,0)=N_D(2s,0)-2$.
It suffices
to rename $s+1=:k_0$ to get (b) (iv).
\end{proof}
Now we are able to prove the following uniqueness result of the complex symplectic structure in the non-unimodular case and in certain unimodular cases:
\begin{corollary}\label{co:uniqueness}
Let $\mfg$ be a $4n$-dimensional almost Abelian Lie algebra admitting a complex symplectic structure. Suppose $\mfg$ is not unimodular or $\mfg$ is unimodular and $\mfg=\mfh\oplus \bR^{4k-1}$ for some $(4(n-k)+1)$-dimensional irreducible almost Abelian Lie algebra $\frh$ and some $k\in \{1,\ldots,n\}$. Then the complex symplectic structure on $\mfg$ is unique up to equivalence.
\end{corollary}
\begin{proof}
Let $(J_1,\omega_1)$ and $(J_2,\omega_2)$ be two complex symplectic structures on $\mfg$. By Proposition \ref{pro:classsimplification}, $(\mfg,J_i,\omega_i)$ is isomorphic to $(\bR^{4n-1}\rtimes_{f_i} \bR,J_0,\omega_0)$ for $i=1,2$ with $f_1$, $f_2$ being as in Proposition \ref{pro:classsimplification}.

Now if $\mfg$ is not unimodular, then $f_1$ and $f_2$ have to be of the first or the third form in Proposition~\ref{pro:classsimplification}. However, since the Lie algebras $\bR^{4n-1}\rtimes_{f_1}\bR$ and $\bR^{4n-1}\rtimes_{f_2}\bR$ are isomorphic as Lie algebras, the real Jordan normal forms of $f_1$ and $f_2$ have to be the same up to scaling. Hence, by (the proof of) Theorem~\ref{th:almostAbelianLAclass}, either both $f_1$ and $f_2$ are of the first form in Proposition~\ref{pro:classsimplification} or both $f_1$ and $f_2$ are of the second form in Proposition~\ref{pro:classsimplification} for the same $m\in \{1,\ldots,n-1\}$ and
the matrices in the left upper corner have to be similar to each other. But these matrices in the left upper corner are symplectically similar by Proposition~\ref{pro:complexJNFsp2mC}~(c) and so the statement follows in the non-unimodular case.

Next, let $\mfg$ be unimodular and $\mfg=\mfh\oplus \bR^{4k-1}$ for some $(4(n-k)+1)$-dimensional irreducible almost Abelian Lie algebra and some $k\in \{1,\ldots,n\}$. Since $\mfh$ is irreducible, the
Jordan normal form of $(f_i)|_{\tilde{\mfu}}$ cannot have any Jordan blocks of size $1$ with zero on the diagonal, where $\tilde{\mfu}$ is an Abelian ideal of codimension one in $\mfh$ and so $\mfu=\tilde{\mfu}\oplus \bR^{4k-1}$. Consequently,
$N_{f_i}(1,0)=4k-1\equiv 3$ (mod 4) and so Theorem \ref{th:almostAbelianLAclass} yields that $f_i$ has to be as in case (II) (i) in that theorem, i.e. $f_i$ is of the second form in Proposition \ref{pro:classsimplification} with $b=c=0$. Thus, here again the assertion follows from Proposition \ref{pro:complexJNFsp2mC} (c) applied to the matrix in $\mathfrak{sp}(2n-2,\bC)$ in the left upper corner of $f_1$ and $f_2$.
\end{proof}

\subsection{Explicit examples of complex symplectic almost Abelian solvmanifolds}\label{subsec:solvmanifold}

We now use the previous results to construct explicit examples of complex symplectic structures on suitable almost Abelian solvmanifolds in any dimension. First of all notice that, in the notations of Theorem \ref{th:complexsymplectic}, an almost Abelian Lie algebra $\bR^{4n-1}\rtimes_f \bR$ admitting a complex symplectic structure is unimodular if and only if $a=0$. By \cite{bock} an almost Abelian Lie group $\bR^{4n-1}\rtimes_\Phi\bR$ admits a lattice if and only if $\exists~t_0\neq 0$ such that $\Phi(t_0)$ is similar to an integer matrix. In this case, a lattice is given by $\Gamma=P^{-1}\mathbb{Z}^{4n-1}\rtimes t_0\mathbb{Z}$, where $P\Phi(t_0)P^{-1}$ is an integer matrix.

We provide examples of complex symplectic almost Abelian solvmanifolds in any dimension. Let $\mathfrak{g}=\bR^{4n-1}\rtimes_f \bR$ be a $4n$-dimensional unimodular almost Abelian Lie algebra with
\begin{equation*}
f=\left(
\begin{array}{c|c}
A & \\
\hline
& \mathbf 0
\end{array}
\right)
	\end{equation*}
where $A\in \mathfrak{sp}(2n-2,\bC)\subseteq \bR^{4(n-1)\times 4(n-1)}$ is the diagonal matrix given by
\[
A\!=\!\text{diag}\left( \frac{1}{2m},\frac{1}{2m},\frac{3}{2m},\frac{3}{2m},\ldots,\frac{2m-1}{2m},\frac{2m-1}{2m},
-\frac{1}{2m}, -\frac{1}{2m},-\frac{3}{2m},-\frac{3}{2m} \ldots,-\frac{2m-1}{2m},-\frac{2m-1}{2m}\right)
\]
with $m\coloneqq n-1$. Then $\mfg$ is as in Theorem \ref{th:almostAbelianLAclass} case (b) (i) and so it admits a complex symplectic structure, which is unique by Corollary \ref{co:uniqueness}. We denote by $G=\bR^{4n-1}\rtimes_\Phi\bR$ the associated simply connected Lie group, where
\begin{equation*}
\Phi(t)=\left(
\begin{array}{c|c}
e^{tA} & \\
\hline
& \mathbf 1
\end{array}
\right)
	\end{equation*}
and
\[
e^{tA}=\text{diag}\,\left( e^\frac{t}{2m},e^\frac{t}{2m},\ldots,e^\frac{(2m-1)t}{2m},
e^\frac{(2m-1)t}{2m},
e^{-\frac{t}{2m}}, e^{-\frac{t}{2m}}, \ldots,e^{-\frac{(2m-1)t}{2m}},e^{-\frac{(2m-1)t}{2m}}\right)\,.
\]
Hence, the characteristic polynomial of $\Phi(t)$ is
\[
P_{\Phi(t)}(x)=(x-1)^3(x-\rho^2)^2(x-\rho^{-2})^2\ldots
(x-\rho^{4m-2})^2(x-\rho^{-(4m-2)})^2,
\]
where $\rho=e^{\frac{t}{4m}}$. Now we argue as in \cite{andrada-origlia}. We set, for $\ell\in\mathbb{N}$, $\ell>2$, $\rho_\ell=\exp\left(\frac{t_\ell}{4m}\right)$, where
\[
t_\ell\coloneqq 2\,m\,\text{arccosh}\left(\frac{\ell}{2}\right)\neq 0.
\]
Then $\rho_\ell^2+\rho_\ell^{-2}=\ell$. We set $a_k=\rho_\ell^{2k}+\rho_\ell^{-2k}$ for $k\geq 0$. Notice that $a_{k+1}=\ell\,a_k-a_{k-1}$ and so $a_k\in\mathbb{Z}$, for every $k$.
Therefore,
\[
(x-\rho_\ell^{2k})(x-\rho_\ell^{-2k})=x^2-a_kx+1
\]
is a polynomial with integer coefficients for every $k$ and so we can write
\[
P_{\Phi(t)}(x)=(x-1)\,q(x)^2,
\]
with
\[
q(x)=x^{2m+1}+b_{2m}x^{2m}+b_{2m-1}x^{2m-1}+\ldots+b_1x-1
\]
for certain $b_1,\ldots,b_{2m}\in \mathbb{Z}$. Thus, $q$
is a polynomial with integer coefficients and distinct roots
and so the corresponding part of $\Phi(t_\ell)$ can be conjugated
to the companion integer matrix
\begin{equation*}
B_q\coloneqq \begin{pmatrix}
 1 & 0 & 0  & 0 & 0 & 1\\
 0 & 0 & 0 & 0 & 0 & 1 \\
0 & 1 & 0 & 0 & 0 & -b_1 \\
 \vdots & \vdots & \ddots & \vdots & \vdots & \vdots \\
 0 & 0 & 0 & 1 & 0 & -b_{2m-1} \\
 0 & 0 & 0 & 0 & 1 & -b_{2m} \\
\end{pmatrix}\,.
\end{equation*}
Thus, $\Phi(t_\ell)$ is conjugated to the integer block-diagonal matrix
\[
B_\ell\coloneqq \diag(1,B_q,B_q)
\]

As a consequence, $G$ admits a lattice
\[
\Gamma_\ell=P_\ell^{-1}\mathbb{Z}^{4n-1}\rtimes_{\Phi} t_\ell\mathbb{Z}
\]
for every $\ell>2$, where $P_\ell\Phi(t_\ell)P_\ell^{-1}=B_\ell$. In particular, $\Gamma_\ell\backslash G$ is a $4n$-dimensional solvmanifold admitting a unique invariant complex symplectic structure. Since $\frg$ is completely solvable, $\Gamma_\ell\backslash G$ admits no K\"ahler metric, see \cite{Hasegawa2}.
\begin{remark}
It can be shown as in \cite[Proposition 4.11]{andrada-origlia} that the solvmanifolds $\Gamma_\ell\backslash G$, with $\ell>2$, are pairwise non homeomorphic.
\end{remark}

\section{Complex symplectic cotangent extension}\label{sec:cotangent}

We adapt to complex symplectic structures the construction of symplectic structures on cotangent extensions proposed in \cite{Ovando}; this will allow us to provide further examples of complex symplectic manifolds. A related construction has been studied in \cite{CPO}. 

Let $\frh$ be a real Lie algebra endowed with a complex structure $J$. On the vector space $\frh^*\oplus\frh$, where $\frh^*=\textrm{Hom}(\frh,\bR)$, we define a skew-symmetric 2-form $\bO$ and an almost complex structure $\bJ$ by
\begin{equation}\label{Om_J}
\bO\big((\varphi,X),(\psi,Y)\big)\coloneqq \varphi(Y)-\psi(X)\,, \qquad \bJ(\varphi,X)\coloneqq (J^*\varphi,JX)\,,
\end{equation}
where $J^*\colon\frh^*\to\frh^*$ is given by~\eqref{J-on-forms}.
It is immediate to see that $\bJ$ is symmetric with respect to $\bO$:
\begin{align*}
\bO\big(\bJ(\varphi,X),(\psi,Y)\big)&=\bO\big((J^*\varphi,JX),(\psi,Y)\big)=(J^*\varphi)(Y)-\psi(JX)=\varphi(JY)-(J^*\psi)(X)\\
&=\bO\big((\varphi,X),(J^*\psi,JY)\big)=\bO\big((\varphi,X),\bJ(\psi,Y)\big)\,.
\end{align*}

Moreover, $\frh^*\subset\frh^*\oplus\frh$ is a complex subspace, meaning that $\bJ\frh^*=\frh^*$, and a Lagrangian subspace, meaning that $\bO\big|_{\frh^*}=0$.

The {\em complex symplectic cotangent extension problem} consists in  finding a Lie algebra structure on $\frh^*\oplus\frh$ such that
\begin{enumerate}
\item[$\mathbf{1.}$] \ $0\to\frh^*\to\frh^*\oplus\frh\to\frh\to 0$ is an exact sequence of Lie algebras, $\frh^*$ being endowed with the structure of an Abelian Lie algebra;
\item[$\mathbf{2.}$] \ $d\bO=0$;
\item[$\mathbf{3.}$] \ $N_\bJ=0$.
\end{enumerate}

A complex symplectic Lie algebra $(\frg,J,\omega)$ is a {\em solution} of the complex symplectic cotangent extension problem if it is isomorphic, as a complex symplectic Lie algebra, to a Lie algebra of the form $(\frh^*\oplus\frh,\bJ,\bO)$.

\medskip

Under the first condition, the most general skew-symmetric bilinear map  
$[\cdot,\cdot]\colon(\frh^*\oplus\frh)\times(\frh^*\oplus\frh)\to\frh^*\oplus\frh$ is determined by a linear map $\rho\colon\frh\to\mathrm{End}(\frh^*)$ and an element $\alpha\in C^2(\frh,\frh^*)=\Lambda^2\frh^*\otimes\frh^*$. Given such $\rho$ and $\alpha$, the bilinear map is
\begin{equation}\label{bracket}
[(\varphi,X),(\psi,Y)] \coloneqq \big(\rho(X)(\psi)-\rho(Y)(\varphi)+\alpha(X,Y),\, [X,Y]_\frh\big)\,;
\end{equation}
here $\varphi,\psi\in\frh^*$ and $X,Y\in \frh$. Where possible, we have shortened the notation $(\varphi,0)$ to $\varphi$, and so on.\\

We now prove the following:
\begin{theorem}\label{thm:cot_extension}
Let $(\frh,J)$ be a Lie algebra endowed with a complex structure. On $\frh^*\oplus\frh$ define a 2-form $\bO$ and an almost complex structure $\bJ$ as in \eqref{Om_J}, and a skew-symmetric bilinear map as in \eqref{bracket}. Then $(\frh^*\oplus\frh,\bJ,\bO)$ is a complex symplectic Lie algebra with Abelian ideal $\frh^*\subset\frh^*\oplus\frh$ if and only if 
\begin{enumerate}[1.]
\item $\alpha\in Z^2(\frh,\frh^*)$;
\item $\rho\colon\frh\to\mathrm{End}(\frh^*)$ is a Lie algebra morphism;
\item $\alpha(X,Y)(Z)+\alpha(Y,Z)(X)+\alpha(Z,X)(Y)=0$;
\item $\alpha(X,Y)-\alpha(JX,JY)=-J^*\big(\alpha(JX,Y)+\alpha(X,JY)\big)$;
\item $\rho(X)(\varphi)(Y)-\rho(Y)(\varphi)(X)+\varphi\big([X,Y]_\frh\big)=0$;
\item $\rho(X)(\varphi)-\rho(JX)(J^*\varphi)=-J^*\big(\rho(X)(J^*\varphi)+\rho(JX)(\varphi)\big)$,
\end{enumerate}
for every $X,Y,Z\in \frh$ and every $\varphi\in\frh^*$.
\end{theorem}

\begin{proof}

First of all, from \cite{Ovando}, we have that the bilinear map $[\cdot,\cdot]$ defined in \eqref{bracket} is a Lie bracket on $\frh^*\oplus\frh$ if and only if $\rho$ is a Lie algebra morphism and $\alpha\in Z^2(\frh,\frh^*)$ is a $2$-cocycle; here the $\frh$-module structure of $\frh^*$ is given by $\rho$. These are conditions $1.$ and $2.$ in the statement. 
Moreover, in \cite{Ovando} it is proved that
$d\bO=0$ if and only if
\begin{align}
\label{Bianchi} \alpha(X,Y)(Z)+\alpha(Y,Z)(X)+\alpha(Z,X)(Y)&=0,\\
\label{rho1}\rho(X)(\varphi)(Y)-\rho(Y)(\varphi)(X)+\varphi\big([X,Y]_\frh\big)&=0\,,
\end{align}
for every $X,Y,Z\in\frh$ and $\varphi\in\frh^*$. These are conditions $3.$ and $5.$ of the theorem; \eqref{Bianchi} is known as {\em Bianchi identity}. We are left with the integrability of $\bJ$. In particular, we prove that $N_\bJ=0$ if and only if
\begin{align}
\label{complex:2}
\alpha(X,Y)-\alpha(JX,JY)&=-J^*\big(\alpha(JX,Y)+\alpha(X,JY)\big)\,,\\
\label{complex:1} 
\rho(X)(\varphi)-\rho(JX)(J^*\varphi)&=-J^*\big(\rho(X)(J^*\varphi)+\rho(JX)(\varphi)\big)
\end{align}
for every $X,Y\in\frh$ and $\varphi\in\frh^*$, giving conditions $4.$ and $6.$ of the theorem. Since we require $\frh^*\subset\frh^*\oplus\frh$ to be an Abelian ideal, we simply need to check the vanishing of $N_\bJ(X,Y)$ and of $N_\bJ(\varphi,X)$, for all $\varphi\in\frh^*$ and $X,Y\in\frh$. For the first case, we have
\begin{align*}
N_\bJ(X,Y)&=[X,Y]+\bJ[JX,Y]+\bJ[X,JY]-[JX,JY]\\
&=\big(\alpha(X,Y),[X,Y]_\frh\big)+\bJ\big(\alpha(JX,Y),[JX,Y]_\frh\big)+\bJ\big(\alpha(X,JY),[X,JY]_\frh\big)\\
&-\big(\alpha(JX,JY),[JX,JY]_\frh\big)\\
&=\big(\alpha(X,Y),[X,Y]_\frh\big)+\big(J^*(\alpha(JX,Y)),J[JX,Y]_\frh\big)+\big(J^*(\alpha(X,JY)),J[X,JY]_\frh\big)\\
&-\big(\alpha(JX,JY),[JX,JY]_\frh\big)\\
&=\alpha(X,Y)+J^*\big(\alpha(JX,Y)\big)+J^*\big(\alpha(X,JY)\big)-\alpha(JX,JY)\,.
\end{align*}
For the second one, we compute
\begin{align*}
N_\bJ(\varphi,X)&=[\varphi,X]+\bJ[J^*\varphi,X]+\bJ[\varphi,JX]-[J^*\varphi,JX]\\
&=-\rho(X)(\varphi)-J^*\big(\rho(X)(J^*\varphi)\big)-J^*\big(\rho(JX)(\varphi)\big)+\rho(JX)(J^*\varphi)\\
&=-\rho(X)(\varphi)-J^*\big(\rho(X)(J^*\varphi)+\rho(JX)(\varphi)\big)+\rho(JX)(J^*\varphi)\,.
\end{align*}
This concludes the proof.
\end{proof}

\begin{remark}\label{rem:1}
We make a few observations on the conditions that appear in Theorem \ref{thm:cot_extension}.
\begin{enumerate}
\item Condition 3. is equivalent to $\alpha\in \Lambda^2 \mfh^*\otimes \mfh^*$ being in the summand $V$ in the irreducible $\GL(\mfh)$ decomposition
\begin{equation*}
\Lambda^2 \mfh^*\otimes \mfh^*=\Lambda^3 \mfh^*\oplus V
\end{equation*}
 with $V=\mathbb{S}_{(2,1)} \mfh^*$; we refer to \cite[Theorem 6.3]{FH} for the notation.
\item If we choose a basis $\{e^1,\ldots,e^{2n}\}$ of $\mfh^*$ with $J^*e^{2j-1}=e^{2j}$ for $j=1,\ldots,n$ and write
\begin{equation*}
\alpha=\sum_{j=1}^{2n} \alpha_j\otimes e^j
\end{equation*}
for uniquely determined $\alpha_1,\ldots,\alpha_{2n}\in \Lambda^2 \mfh^*$, then it is easy to see that condition 4. is equivalent to $\alpha_{2j-1}+i \alpha_{2j}$ having no $(0,2)$-part, which in turn is equivalent to
\begin{equation*}
\alpha_{2j-1}=\sigma_j+\re(\psi_j),\qquad \alpha_{2j}=\tau_j+\im(\psi_j)
\end{equation*}
for real $(1,1)$-forms $\sigma_j,\tau_j\in [\Lambda^{1,1}\mfh^*]$ and complex $(2,0)$-forms $\psi_j\in \Lambda^{2,0}\mfh^*$ for $j=1,\ldots,n$.
\item If $\rho=0$, then condition 5. forces $\mfh$ to be nilpotent and $\frh^*\oplus\frh$ is a particular kind of central extension of the Abelian Lie algebra $\frh$ by the cocycle $\alpha$.
\end{enumerate}
\end{remark}

Notice that a complex symplectic Lie algebra $(\frg,J,\omega)$ of real dimension $4n$ with an Abelian, $J$-invariant $2n$-dimensional ideal $\frh^*$ is a solution of the complex symplectic cotangent extension problem if and only if conditions 1. $-$ 6. hold.

\begin{example}\label{ex:cot_extension}
The 6-dimensional nilpotent Lie algebra $\frh:=\frh_7=(0,0,0,12,13,23)$ has only one complex structure $J$ up to isomorphisms (see \cite{COUV}), namely
\[
J=-e^1\otimes e_2+e^2\otimes e_1+e^3\otimes e_4-e^4\otimes e_3-e^5\otimes e_6+e^6\otimes e_5\,.
%
\]
If we want $\frg=\frh^*\oplus\frh$ to be nilpotent, we need the matrices $\rho$ to be nilpotent, hence strictly upper triangular with respect to some basis of $\frh^*$; in our setting, the natural basis to work with is the basis dual to $\{e_1,\ldots,e_6\}$. A direct computation shows that the most general such $\rho$ satisfying \eqref{rho1} and \eqref{complex:1} consists of the following matrices, written in the coframe $\{e^1,\ldots,e^6\}$ of $\frh^*$:
\[
\rho(e_1)=\begin{pmatrix}
0 & 0 & \rho^1_{13} & \rho^1_{14} & \rho^1_{15} & \rho^1_{16}\\
0 & 0 & \rho^1_{23} & \rho^1_{24} & \rho^1_{25} & \rho^1_{26}\\
0 & 0 & 0 & 0 & \rho^1_{35} & \rho^1_{36}\\
0 & 0 & 0 & 0 & \rho^1_{45} & \rho^1_{46}\\
0 & 0 & 0 & 0 & 0 & 0\\
0 & 0 & 0 & 0 & 0 & 0
\end{pmatrix}\,, \quad
\rho(e_2)=\begin{pmatrix}
0 & 0 & \rho^1_{23} & \rho^1_{24}-1 & \rho^1_{25} & \rho^1_{26}\\
0 & 0 & -\rho^1_{13}-1 & -\rho^1_{14} & -\rho^1_{15} & -\rho^1_{16}\\
0 & 0 & 0 & 0 & -\rho^1_{45} & -\rho^1_{46}\\
0 & 0 & 0 & 0 & \rho^1_{35} & \rho^1_{36}\\
0 & 0 & 0 & 0 & 0 & 0\\
0 & 0 & 0 & 0 & 0 & 0
\end{pmatrix}\,,
\]
\[
\rho(e_3)=\begin{pmatrix}
0 & 0 & 0 & 0 & \rho^1_{35}-1 & \rho^1_{36}\\
0 & 0 & 0 & 0 & -\rho^1_{45} & -\rho^1_{46}-1\\
0 & 0 & 0 & 0 & \rho^3_{35} & \rho^3_{36}\\
0 & 0 & 0 & 0 &  \rho^3_{45} & \rho^3_{46} \\
0 & 0 & 0 & 0 & 0 & 0\\
0 & 0 & 0 & 0 & 0 & 0
\end{pmatrix}\,, \quad
\rho(e_4)=\begin{pmatrix}
0 & 0 & 0 & 0 & \rho^1_{45} & \rho^1_{46}\\
0 & 0 & 0 & 0 & \rho^1_{35} & \rho^1_{36}\\
0 & 0 & 0 & 0 & \rho^3_{45} & \rho^3_{46}\\
0 & 0 & 0 & 0 & -\rho^3_{35}& -\rho^3_{36}\\
0 & 0 & 0 & 0 & 0 & 0\\
0 & 0 & 0 & 0 & 0 & 0
\end{pmatrix}\,,
\]
together with $\rho(e_5)=\rho(e_6)=\mathbf 0$. Imposing that $\rho\colon\frh\to\mathrm{End}(\frh^*)$ is a Lie algebra morphism yields the linear equations $\rho^3_{35}=\rho^3_{36}=\rho^3_{45}=\rho^3_{46}=0$ and the following four quadratic equations:
\begin{equation*}
	\begin{split}
	&(\rho^1_{13}+\rho^1_{24}-2)\,\rho^1_{45}+(\rho^1_{23}-\rho^1_{14})\,\rho^1_{35}=0,\\
	&(\rho^1_{23}-\rho^1_{14})\,\rho^1_{45}-(\rho^1_{13}+\rho^1_{24}+2)\,\rho^1_{35}=0,
	\end{split}
\qquad
	\begin{split}
	&(\rho^1_{13}+\rho^1_{24}-2)\,\rho^1_{46}+(\rho^1_{23}-\rho^1_{14})\,\rho^1_{36}=0,\\
	&(\rho^1_{23}-\rho^1_{14})\,\rho^1_{46}-(\rho^1_{13}+\rho^1_{24}+2)\,\rho^1_{36}=0.
	\end{split}
\end{equation*}
We find four sets of solutions:
\[
\left\{
\begin{array}{ccl}
\rho^1_{23} & = & \rho^1_{14}\\
\rho^1_{24} & = & -2-\rho^1_{13}\\
\rho^1_{45} & = & 0\\
\rho^1_{46} & = & 0
\end{array}
\right.\quad
\left\{
\begin{array}{ccl}
\rho^1_{23} & = & \rho^1_{14}\\
\rho^1_{24} & = & 2-\rho^1_{13}\\
\rho^1_{35} & = & 0\\
\rho^1_{36} & = & 0\\
\end{array}
\right. \quad
\left\{
\begin{array}{ccl}
\rho^1_{35} & = & 0\\
\rho^1_{36} & = & 0\\
\rho^1_{45} & = & 0\\
\rho^1_{46} & = & 0
\end{array}
\right.
\quad \mathrm{and} \quad
\left\{
\begin{array}{ccl}
\rho^1_{23} & = & \rho^1_{14}+2\sin t\\
\rho^1_{24} & = & -\rho^1_{13}+2\cos t\\
\rho^1_{45} & = & \frac{\cos t+1}{\sin t}\rho^1_{35}\\
\rho^1_{46} & = & \frac{\cos t+1}{\sin t}\rho^1_{36}
\end{array}
\right.\,,
\]
for $t\in(0,2\pi)$, $t\neq\pi$; the parameters which do not appear are free. Thus we find many nilpotent solutions of the cotangent extension problem. In this example we work with $\alpha=0$; the presence of $\alpha$ does not alter the nilpotency condition, as it simply gives a particular kind of central extension.
\end{example}

Suppose $\frh$ is a Lie algebra endowed with a complex structure $J$, and let $\bJ$ be the complex structure on the cotangent extension $\frh^*\oplus\frh$ given in \eqref{Om_J}; we obtain necessary and sufficient conditions for $\bJ$ to be Abelian or parallelizable.


\begin{proposition}\label{prop:ab_par}
Let $\frh$ be a Lie algebra endowed with a complex structure $J$ and let $\bJ$ be the complex structure given by \eqref{Om_J} on $\frh^*\oplus\frh$, with Lie algebra structure determined by $\rho$ and $\alpha$ as above. Then
\begin{itemize}
\item $\bJ$ is Abelian if and only if 
\begin{itemize}
\item $J$ is Abelian;
\item $\alpha\in[\Lambda^{1,1}\frh^*]\otimes\frh^*$;
\item $\rho\colon\frh\to\mathrm{End}(\frh^*)$ is anti-holomorphic;
\end{itemize}
\item $\bJ$ is parallelizable if and only if
\begin{itemize}
\item $J$ is parallelizable;
\item $J^*(\alpha(X,Y))=\alpha(JX,Y)$, $\forall~X,Y\in\frh$;
\item $\rho\colon\frh\to\mathrm{End}^{J^*}(\frh^*)$;
\item $J^*\circ\rho(X)=\rho(JX)$ for all $X\in\frh$.
\end{itemize}
\end{itemize}
\end{proposition}

An (almost) complex structure $\hat{J}$ on $\mathrm{End}(\frh^*)$ is defined by $\hat{J}(f)(\varphi)\coloneqq f(J^*\varphi)$. That $\rho\colon\frh\to\mathrm{End}(\frh^*)$ is anti-holomorphic means then that $\rho\circ J=-\hat{J}\circ\rho$. Also, $\mathrm{End}^{J^*}(\frh^*)\coloneqq\{f\in\mathrm{End}(\frh^*) \mid f\circ J^*=J^*\circ f\}$.

%
%
%
%

\begin{proof}
Let us consider the Abelian case. Computing as in the final part of the proof of Theorem \ref{thm:cot_extension}, we see that, for $X,Y\in\frh$,
\[
[X,Y]-[JX,JY]=(\alpha(X,Y)-\alpha(JX,JY),[X,Y]_\frh-[JX,JY]_\frh)\,,
\]
which vanishes if and only if $J$ is Abelian and $\alpha(X,Y)=\alpha(JX,JY)$, that is, $\alpha\in[\Lambda^{1,1}\frh^*]\otimes\frh^*$. For $\varphi\in\frh^*$ and $X\in\frh$ we have
\[
[\varphi,X]-[J^*\varphi,JX]=-\rho(X)(\varphi)+\rho(JX)(J^*\varphi)\,,
\]
which vanishes if and only if $\rho(JX)(J^*\varphi)=\rho(X)(\varphi)$; this is equivalent to $\hat{J}(\rho(JX))(\varphi)=\rho(X)(\varphi)$ for all $\varphi\in\frh^*$, and to $\rho(JX)=-\hat{J}(\rho(X))$, which in turns means $\rho\colon\frh\to\mathrm{End}(\frh^*)$ is anti-holomorphic.

As for the parallelizable case, we compute
\[
\bJ[X,Y]-[JX,Y]=(J^*(\alpha(X,Y))-\alpha(JX,Y),J[X,Y]_\frh-[JX,Y]_\frh)\,,
\]
which vanishes if and only if $J$ is parallelizable and $J^*(\alpha(X,Y))=\alpha(JX,Y)$, for all $X,Y\in\frh$. On the other hand,
\[
\bJ[\varphi,X]-[J^*\varphi,X]=-J^*(\rho(X)(\varphi))+\rho(X)(J^*\varphi)\,,
\]
which vanishes if and only if $\rho(X)\in\mathrm{End}^{J^*}(\frh^*)$. Finally,
\[
\bJ[X,\varphi]-[J X,\varphi]=J^*(\rho(X)(\varphi))-\rho(JX)(\varphi)\,,
\]
which vanishes if and only if $J^*\circ\rho(X)=\rho(JX)$ for all $X\in\frh$.
\end{proof}

\begin{remark}
Condition \eqref{complex:2} above holds under either of the assumptions on $\alpha$ made in Proposition \ref{prop:ab_par}. Similarly, condition \eqref{complex:1} holds under each of the assumptions on $\rho$ made in Proposition \ref{prop:ab_par}.
\end{remark}

We aim now at characterizing complex symplectic Lie algebras $(\frg,J,\omega)$ which are solutions of the cotangent extension problem. A necessary condition for this is that $\frg$ admits a Lagrangian (hence Abelian, by \cite[Lemma 4.1]{CPO}) $J$-invariant ideal. We will show that this condition is also sufficient. Suppose $(\frg,J,\Omega)$ is a complex symplectic Lie algebra and assume that $\frj\subset\frg$ is a Lagrangian ideal; then one can invoke \cite{Weinstein} to conclude that $\frj$ has a Lagrangian complement in $\frg$. If, in addition, $\frj$ is $J$-invariant, then it is Lagrangian with respect to the complex symplectic 2-form $\omega_\bC$. We show that this is the case: one can take $\frj$ to be $J$-invariant.

\begin{lemma}\label{lag:complement}
Let $V$ be a $4n$-dimensional vector space endowed with a symplectic structure $\omega$ and an almost complex structure $J$ such that $J$ is either symmetric or skew-symmetric with respect to $\omega$, i.e. either $\omega(Ju,v)=\omega(u,Jv)$, or $\omega(Ju,v)=-\omega(u,Jv)$ $\forall\,u,v\in V$. Assume $L\subset V$ is a Lagrangian, $J$-invariant subspace. Then there exists a Lagrangian, $J$-invariant subspace $L'$ such that $V=L\oplus L'$.
\end{lemma}

\begin{proof}
We adapt ideas from \cite[\S 2.2]{GS}. If $\tilde{g}$ is a scalar product on $V$, then $g(v,w)\coloneqq \tilde{g}(v,w)+\tilde{g}(Jv,Jw)$ is a $J$-invariant scalar product. $L^\perp$, the $g$-orthogonal of $L$, is a $J$-invariant subspace of $V$, of dimension $2n$. Any $J$-invariant $2n$-dimensional subspace $W\subset V$ such that $W\cap L=\{0\}$ is the graph of a linear map $f\colon L^\perp\to L$ such that $f\circ J=J\circ f$; we denote by $\mathrm{Hom}_J(L^\perp,L)$ the space of such maps. In other words, by writing $V=L^\perp\oplus L$, the elements of $W$ are of the form $u+f(u)$, for $u\in L^\perp$. We have
\begin{equation}\label{eq:242}
\omega(u+f(u),v+f(v))=\omega(u,v)+\omega(u,f(v))+\omega(f(u),v)\,,
\end{equation}
because $\omega(f(u),f(v))=0$, since $f(u),f(v)\in L$. Now to any $f\in\mathrm{Hom}(L^\perp,L)$ we can associate a bilinear form $b_f$ on $L^\perp$ via the formula 
\begin{equation}\label{eq:243}
b_f(u,v)\coloneqq\omega(f(u),v)\,.
\end{equation}
We have
\[
\mathrm{Hom}(L^\perp,L)\cong L\otimes (L^\perp)^*\cong (L^\perp)^*\otimes (L^\perp)^*\,,
\]
where the identification between $L$ and $(L^\perp)^*$ in the last step is given by $\omega$. Hence, the map $\mathrm{Hom}(L^\perp,L)\to\mathrm{Bil}(L^\perp)$, $f\mapsto b_f$, given by \eqref{eq:243} is an isomorphism. Now suppose $J$ is symmetric with respect to $\omega$. Then $f\in\mathrm{Hom}_J(L^\perp,L)$ if and only if $J$ is symmetric with respect to $b_f$. Indeed, if $f\in\mathrm{Hom}_J(L^\perp,L)$, then
\[
b_f(Ju,v)=\Omega(f(Ju),v)=\Omega(Jf(u),v)=\Omega(f(u),Jv)=b_f(u,Jv)\,.
\] 
On the other hand, if $J$ is symmetric with respect to $b_f$, then, $\forall\,v\in L^\perp$,
\[
\omega(f(Ju),v)=b_f(Ju,v)=b_f(u,Jv)=\omega(f(u),Jv)=\omega(Jf(u),v)\,,
\]
and $f\in \mathrm{Hom}_J(L^\perp,L)$. We obtain an analogous conclusion when $J$ is skew-symmetric with respect to $\omega$: $f\in\mathrm{Hom}_J(L^\perp,L)$ if and only if $J$ is skew-symmetric with respect to $b_f$. Going back to \eqref{eq:242}, $W$ is Lagrangian if and only if
\[
b_f(u,v)-b_f(v,u)=-\omega(u,v)\,.
\]

Choosing $f\in \mathrm{Hom}(L^{\perp},L)$ with $b_f=-\tfrac{1}{2}\omega$ solves this equation; we only have to make sure that $f$ commutes with $J$, that is, $f\in\mathrm{Hom}_J(L^{\perp},L)$. But if $J$ is symmetric/skew-symmetric with respect to $\omega$, the same holds for $b_f$, and so by what we showed above, we have $f\in \mathrm{Hom}_J(L^{\perp},L)$.


\end{proof}

We have the following result:
\begin{theorem}\label{theo:cot_extension}
A complex symplectic Lie algebra $(\frg,J,\omega)$ is a solution of the complex symplectic cotangent extension problem if and only if it admits a $J$-invariant Lagrangian ideal.
\end{theorem}

\begin{proof}
The necessity has already been discussed. Consider the short exact sequence of Lie algebras
\begin{equation}\label{eq:short-exact-seq}
0\to\frj\to\frg\to\frh\to 0\,,
\end{equation}
where $\frh\coloneqq \frg/\frj$ is the quotient Lie algebra; since $\frj$ is $J$-invariant, $J$ induces a complex structure on $\frh$, which we denote again by $J$. Since $\frj$ is Lagrangian, it is Abelian; moreover, since $\omega$ is non-degenerate, the linear map $\sigma\colon\frj\to\frh^*$, $u\mapsto \imath_u\omega$ is injective, hence an isomorphism by dimension reasons; moreover, it satisfies $\sigma\circ J=J^*\circ\sigma$. Since $J$ is symmetric with respect to $\omega$, we can use Lemma \ref{lag:complement} to produce a Lagrangian, $J$-invariant complement $\frl$, so that $\frg=\frj\oplus\frl$. The identification $\tau\colon\frl\to\frh$ commutes with $J$. We use $\sigma$ and $\tau$ to endow $\frh^*\oplus\frh$ with a Lie algebra structure so that $\sigma\oplus\tau\colon\frg=\frj\oplus\frl\to\frh^*\oplus\frh$ is a Lie algebra isomorphism which commutes with the respective complex structures. Clearly $\frh^*$ sits
in $\frh^*\oplus\frh$ as an Abelian ideal. If we endow $\frh^*\oplus\frh$ with the canonical complex symplectic structure $(\bJ,\bO)$, it is easy to see that $\sigma\oplus\tau\colon(\frg,J,\omega)\to(\frh^*\oplus\frh,\bJ,\bO)$ is an isomorphism of complex symplectic Lie algebras.
\end{proof}

\subsection{Explicit solutions to the cotangent extension problem}

In this section we work out some explicit solutions to the complex symplectic cotangent extension problem.

\subsubsection{The case $\rho=0$}\label{subsec:rho=0}

In this case, $\mfh$ has to be Abelian as we noticed in Remark \ref{rem:1}. As a consequence of the short exact sequence~\eqref{eq:short-exact-seq} in the proof of Theorem~\ref{theo:cot_extension}, one has
\begin{proposition}
 Let $(\mfg,J,\omega)$ be a complex symplectic Lie algebra admitting a $J$-invariant Lagrangian ideal $\mfj$. If $\mfj$ is central, then $\mfg/\mfj$ is Abelian.
\end{proposition}

Now note that $\rho=0$ and $\mfh$ being Abelian imply that conditions 2., 5. and 6. in Theorem \ref{thm:cot_extension} are automatically satisfied and that condition 1. is valid for any $\alpha\in \Lambda^2 \mfh^*\otimes \mfh^*$. Hence, $\alpha$ only needs to satisfy conditions 3. and 4. in Theorem \ref{thm:cot_extension}  (or, equivalently, conditions (1) and (2) in Remark \ref{rem:1}) in order to define a solution to the complex symplectic cotangent extension problem.
Since we are dealing with a central extension in this case, $\mfg$ is two-step nilpotent. Notice that not all solutions with $\mfg$ two-step nilpotent necessarily have $\rho=0$, as one sees, for instance, using the last set of solutions for $\rho$ in Example \ref{ex:cot_extension}. Moreover, note that the complex structure $\mathbf{J}$ on $\mfg$ is Abelian if and only if $\alpha\in[\Lambda^{1,1}\mfh^*]\otimes \mfh^*$. 
\begin{example} Let us now give some examples.
\begin{itemize}
	\item[(a)]
	The simplest idea is to take $\alpha_{2j-1}=e^{2j-1}\wedge e^{2j}$, $\alpha_{2j}=0$ for $j=1,\ldots,n$. Then $\mathbf{J}$ is Abelian and $\mfg$ is isomorphic as a Lie algebra to $(\mfh_3\oplus \bR)^n$.
	
	Note that in the case $n=1$, up to isomorphism, the data for $\alpha$ as above is the only possible solution of the complex symplectic cotangent extension problem with $\rho=0$ and $\mfg$ not being Abelian.
	\item[(b)]
	We consider now the case $n=2$ and determine first all possible choices of $\alpha$. For this, we set
	\begin{equation*}
	\omega_1:=e^{12},\ \omega_2:=e^{34},\ \omega_3:=e^{13}+e^{24},\ \omega_4:=e^{14}-e^{23},\ \sigma_1:=e^{13}-e^{24},\ \sigma_2=e^{14}+e^{23}.
	\end{equation*}
	Condition (3) from Remark \ref{rem:1} is satisfied if and only if
	\begin{equation*}
	\begin{split}
	\alpha_1&=\sum_{i=1}^4 a_i\, \omega_i+ a_5\, \sigma_1+a_6\, \sigma_2,\quad 	\alpha_2=\sum_{i=1}^4 b_i\, \omega_i+ a_6\, \sigma_1-a_5\,  \sigma_2,\\
		\alpha_3&=\sum_{i=1}^4 c_i\, \omega_i+ c_5\, \sigma_1+c_6\, \sigma_2,\quad 	\alpha_4=\sum_{i=1}^4 d_i\, \omega_i+ c_6\, \sigma_1-c_5\,  \sigma_2.
	\end{split}
	\end{equation*}
	for certain $a_1,\ldots,a_6,b_1,\ldots,b_4,c_1,\ldots,c_6,d_1,\ldots,d_4\in \bR$. Then, condition 2. from Remark \ref{rem:1} is fulfilled if and only if
	\begin{equation*}
	c_1=b_3+a_4,\quad d_1=b_4-a_3,\quad a_2=c_4-d_3,\quad b_2=c_3+d_4.
	\end{equation*}
	We discuss some special cases:
\begin{itemize}
	\item[(i)]
	Choosing $\alpha_1=\sigma_1$, $\alpha_2=-\sigma_4$, $\alpha_3=\alpha_4=0$, we get a solution of the complex symplectic cotangent extension problem on $\mfg=\mfh_3(\bC)\oplus \bR^2$ with non-Abelian complex structure.
	\item[(ii)]
	Choosing $\alpha_1=\sigma_1$, $\alpha_2=-\sigma_4$, $\alpha_3=\omega_2$, $\alpha_4=0$, we get a solution of the complex symplectic cotangent extension problem on $\mfg=(37B_1)\oplus \bR$ with non-Abelian complex structure; here $(37B_1)$ is the notation used by Gong in \cite{Gong}.
	\item[(iii)]
	Choosing $\alpha_1=\omega_4$, $\alpha_2=\omega_3$, $\alpha_3=2 \omega_1$ and $\alpha_4=\delta \omega_2$ for some $\delta\in \{0,1\}$, we get a solution of the complex symplectic cotangent extension problem with Abelian complex structure. For $\delta=0$, the underlying Lie algebra $\mfg$ is again $(37B_1)\oplus \bR$, whereas for $\delta=1$, the Lie algebra $\mfg$ is an indecomposable eight-dimensional two-step nilpotent Lie algebra.
\end{itemize}
\end{itemize}
\end{example}

\subsubsection{$\mfh$ Abelian with $\hrho$ of full rank}
Here, $\hrho$ is the linear map $\rho\colon \mfh\rightarrow \End(\mfh^*)$ considered as a map $\hrho:\mfh\otimes \mfh\rightarrow \mfh$ via the relation $\varphi(\hrho(X,Y))=\rho(X)(\varphi)(Y)$ for $X,Y\in \mfh$, $\varphi\in \mfh^*$.

\begin{claim}
$\hrho$ has full rank if and only if for any $\varphi\in \mfh^*\setminus \{0\}$ there exist some $X\in \mfh$ with $\rho(X)(\varphi)\neq 0$.
\end{claim}

Indeed, if $\hrho$ has full rank and $\varphi\in\mfh^*\setminus \{0\}$ is given, choose some $Z\in \mfh$ with $\varphi(Z)\neq 0$. Since $\varphi$ has full rank, there exist $X_i,Y_i\in \mfh$, $i=1,\ldots, N$ with
$\sum_{i=1}^N \hrho(X_i,Y_i)=Z$ and so
\begin{equation*}
0\neq\varphi(Z)=\sum_{i=1}^N \varphi(\hrho(X_i,Y_i))=\sum_{i=1}^N \rho(X_i)(\varphi)(Y_i)
\end{equation*}
and so $\rho(X_i)(\varphi)\neq 0$ for at least one $i\in \{1,\ldots,N\}$.

Conversely, if $\hrho$ does not have full rank, there exists some $\varphi\in \mfh^*\setminus \{0\}$ such that $\varphi(\mathrm{im}(\hrho))=\{0\}$ and so $0=\varphi(\hrho(X,Y))=\rho(X)(\varphi)(Y)$ for all $X,Y\in \mfh$. But then $\rho(X)(\varphi)=0$ for all $X\in \mfh$.
\begin{flushright}
\qedsymbol{}
\end{flushright}

Note that since $\mfh^*$ is an Abelian ideal in $\mfg$, the claim is equivalent to $\mfh^*\cap \mathfrak{z}(\mfg)=\{0\}$. We next use this characterization to prove the following
\begin{theorem}
Let $\mfh$ be the Abelian Lie algebra of dimension $2n$ with standard complex structure $J$ and let $\rho:\mfh\rightarrow \End(\mfh^*)$ be linear such that the associated linear map $\hrho:\mfh\otimes \mfh\rightarrow \mfh$ has full rank. Then $\rho$, together with $\alpha=0$, gives a solution to the complex symplectic cotangent extension problem if and only if $\hrho\in S^2\mfh^*\otimes \mfh$, $\hrho(\hrho(\cdot,\cdot),\cdot)\in S^3\mfh^*\otimes \mfh$ and $\hrho(J\cdot,\cdot)=\hrho(\cdot,J\cdot)$. In this case, there exists $X\in \mfh$ such that $\hrho(X,\cdot)=id_{\mfh}$ and $\hrho(JX,\cdot)=J$.
\end{theorem}

\begin{proof}
For $\alpha$ not necessarily equal to zero, we observe the following.
First, since $\mfh$ is Abelian, condition 5. in Theorem \ref{thm:cot_extension} is equivalent to $\hrho$ being symmetric in its arguments, i.e.~$\hrho\in S^2\mfh^*\otimes \mfh$. Moreover, the condition that $\rho$ is a representation of $\mfh$ on $\mfh^*$ is equivalent to $\hrho(\hrho(\cdot,\cdot),\cdot)\in S^3\mfh^*\otimes \mfh$. 

Since $\mathrm{im}(\hrho)=\mfh$, \cite[Proposition 4.3]{FrSw} yields the existence of some $X\in \mfh$ with $\hrho(X,\cdot)=\id_{\mfh}$, i.e. $\rho(X)=\id_{\mfh^*}$. Using again that $\mfh$ is Abelian, condition 1. in Theorem \ref{thm:cot_extension} implies
 \begin{equation*}
 \rho(X)(\alpha(Y,Z))+\rho(Y)(\alpha(X,Z))+\rho(Z)(\alpha(X,Y))=0,
 \end{equation*}
 i.e.
 \begin{equation}\label{eq:alphaexact}
\alpha(Y,Z)=\rho(Y)(\alpha(X,Z))-\rho(Z)(\alpha(X,Y))=\rho(Y)(\nu(Z))-\rho(Z)(\nu(Y))
 \end{equation}
 for all $Y,Z\in \mfh$, where $\nu\in\frh^*\otimes\frh^*=\Lambda^1(\mfh,\mfh^*)$ is defined by $\nu(Y)\coloneqq\alpha(X,Y)$ for $Y\in \mfh$.
 Equation \eqref{eq:alphaexact} says that $\alpha=d_{\rho}\nu$, i.e.~$\alpha$ is $d_{\rho}$-exact, thus the short exact sequence 
\begin{equation*}
0\rightarrow \mfh^*\rightarrow \mfh^*\oplus \mfh\rightarrow \mfh\rightarrow 0
\end{equation*}
splits, and $\mfg=\mfh^*\rtimes \mfh$. 

Let us now assume that $\alpha=0$. We have to find $\hrho\in S^2\mfh^*\otimes \mfh$ such that $\hrho(\hrho(\cdot,\cdot),\cdot)\in S^3\mfh^*\otimes \mfh$ and condition 6. of Theorem \ref{thm:cot_extension} holds; expressed in terms of $\hrho$, it amounts to
\begin{equation*}
\hrho(Y,Z)+\hrho(JY,JZ)+J(\hrho(Y,JZ)-\hrho(JY,Z))=0
\end{equation*}
for all $Y,Z\in \mfh$. In particular, taking $Z=JY$, we obtain
\begin{equation*}
0=\hrho(Y,JY)-\hrho(JY,Y)+J(-\hrho(Y,Y)-\hrho(JY,JY))=-J(\hrho(Y,Y)+\hrho(JY,JY))\,,
\end{equation*}
i.e. $\hrho(JY,JY)=-\hrho(Y,Y)$ for all $Y\in \mfh$. By polarizing, we obtain 
\begin{equation}\label{eq:101}
\hrho(JY,JZ)=-\hrho(Y,Z)
\end{equation}
for all $Y,Z\in \mfh$. 

Conversely, \eqref{eq:101} implies condition 6. in Theorem \ref{thm:cot_extension}. In particular, $\hrho(JX,Y)=\hrho(X,JY)=JY$, i.e. $\hrho_{JX}\coloneqq \hrho(JX,\cdot)=J$. 

\end{proof}

Finally observe that $\mfh$ being Abelian implies $[\mfg,\mfg]\subseteq \mfh^*$, hence that $\mfg$ is two-step solvable. As a consequence, we thus have obtained

\begin{corollary}\label{th:rhofullrank}
Let $(\mfg,J,\omega)$ be a complex symplectic Lie algebra admitting a $J$-invariant Lagrangian ideal $\mfj$. If $\mfj\cap \mathfrak{z}(\mfg)=\{0\}$ and $\mfg/\mfj$ is Abelian, then $\mfg$ is two-step solvable and there exist an Abelian subalgebra $\mfa$ of $\mfg$ such that $\mfg=\mfj\rtimes \mfa$ as Lie algebras, as well as some $X\in \mfa$ with $\ad(X)|_{\mfj}=\id_{\mfj}$ and $\ad(JX)|_{\mfj}=J|_{\mfj}$.
\end{corollary}

Next, we give a complete classification in low dimensions:
\begin{theorem}
Let $(\mfg,J,\omega)$ be a complex symplectic Lie algebra admitting a $J$-invariant Lagrangian ideal $\mfj$. Furthermore, assume that $\mfj\cap \mathfrak{z}(\mfg)=\{0\}$ and $\frh\coloneqq\mfg/\mfj$ is Abelian. Then:
\begin{itemize}
\item[a)]
If $\dim(\mfg)=4$, then, up to isomorphism, $(J,\omega)$ is the complex symplectic structure on $\mfr_2'$ from \cite{BGGL}.
\item[b)]
If $\dim(\mfg)=8$, then $\mfg$ has a basis $\{f_1,\ldots,f_8\}$ such that the differentials of the dual basis are given by
\begin{equation*}
(0^4,-15+26-37+48,-16-25-38-47,-17+28-\delta(35-46),-18-27-\delta(36+45)),
\end{equation*}
for some $\delta\in \{0,1\}$ and such that
\begin{equation*}
J=\sum_{i=1}^4 \left(f^{2i-1}\otimes f_{2i}-f^{2i}\otimes f_{2i-1}\right),\quad
\omega=f^{15}-f^{26}+f^{37}-f^{48}\,.
\end{equation*}
\end{itemize}
\end{theorem}
\begin{proof}
For a), let $e_1,e_2=Je_1$ be a basis of $\mfh$ such that $\hrho(e_1,\cdot)=\id_{\mfh}$.
Then $(e_1,e_2,e^1,-e^2)$ is a basis of $\mfg$. If we denote this basis by $f_1,\ldots,f_4$, then, up to skew-symmetry, Theorem \ref{th:rhofullrank} shows that the only non-zero Lie brackets are given by
\begin{equation*}
[f_1,f_3]=f_3,\quad [f_1,f_4]=f_4,\quad [f_2,f_3]=f_4,\quad [f_2,f_4]=-f_3,
\end{equation*}
which shows that $\mfg$ is isomorphic as a Lie algebra to $\mfr_2'$. Moreover, we have
\begin{equation*}
J=f^1\otimes f_2-f^2\otimes f_1+f^3\otimes f_4-f^4\otimes f_3,\quad \omega=f^{13}-f^{24}
\end{equation*}
and so $(J,\omega)$ is the complex symplectic structure on $\mfr_2'$ from \cite{BGGL}.

In order to prove b), we choose again $e_1,e_2\coloneqq Je_1$ such that $\hrho(e_1,\cdot)=\id_{\mfh}$ and then $\hrho(e_2,\cdot)=J|_{\mfh}$. We extend $e_1,e_2$ to a basis $e_1,\ldots,e_4$ such that $Je_3=e_4$. We first show that we may do this in such a way that $\hrho(e_j,e_k)\in \spa{e_1,e_2}$ for$j,k\in \{3,4\}$. For this, we consider
\begin{equation*}
e_3':= e_3+ a e_1+ be_2,\quad e_4':=Je_3'
\end{equation*}
for $a,b\in \bR$. Then, it suffices to find $a,b\in \bR$ such that $\hrho(e_3',e_3')\in \spa{e_1,e_2}$; indeed, if this is the case, then $\hrho(e_4',e_4')=-\hrho(e_3',e_3')$ and
\begin{equation*}
\hrho(e_4',e_3')=\hrho(e_3',e_4')=\hrho(e_3',Je_3')=\hrho(e_3',\hrho(e_2,e_3'))=\hrho(e_2,\hrho(e_3',e_3'))=J \hrho(e_3',e_3')\,.
\end{equation*}
We compute
\begin{equation*}
\begin{split}
\hrho(e_3',e_3')&=\hrho(e_3+a e_1+b e_2,e_3+a e_1+ b e_2)\\
&=\hrho(e_3,e_3)+a\hrho(e_3,e_1)+b\hrho(e_3,e_2)+a\hrho(e_1,e_3+a e_1+ b e_2)+b\hrho(e_2,e_3+a e_1+ b e_2)\\
&=\hrho(e_3,e_3)+2 a\hrho(e_1,e_3)+ 2b\hrho(e_2,e_3)+a^2\hrho(e_1,e_1)+2ab\hrho(e_1,e_2)+b^2\hrho(e_2,e_2)\\
&=\hrho(e_3,e_3)+2 a e_3 + 2b e_4+ (a^2-b^2) e_1+ 2ab e_2
\end{split}
\end{equation*}
and see that we may choose $a,b\in \bR$ such that $\hrho(e_3',e_3')\in \spa{e_1,e_2}$. Thus, we may assume that there exist $c,d\in \bR$ with
\begin{equation*}
\hrho(e_3,e_3)=c e_1+ de_2.
\end{equation*}
By the computations above, we then have
\begin{equation*}
\hrho(e_4,e_3)=\hrho(e_3,e_4)= J\hrho(e_3,e_3)=-d e_1+ c e_2\,,\quad \hrho(e_4,e_4)=-\hrho(e_3,e_3)=-c e_1-d e_2.
\end{equation*}
Now we show that either $\hrho(e_j,e_k)=0$ for all $j,k\in \{3,4\}$, or
we may change the basis in such a way that $\hrho(e_3,e_3)=e_1$.

The first case occurs if $(c,d)=(0,0)$. So let us now assume that $(c,d)\neq (0,0)$.
We may assume that $d\neq 0$ since for $d=0$, we already have $\hrho(e_3,e_3)=ce_1$ and so may normalise $e_3$ such that $\hrho(e_3,e_3)=e_1$. For $d\neq 0$, we take $e_3'\coloneqq e_3+\mu e_4$ and compute
\begin{equation*}
\begin{split}
\hrho(e_3',e_3')&=\hrho(e_3+\mu e_4,e_3+\mu e_4)=c e_1+ d e_2+ 2\mu (-d e_1+ c e_2)- \mu^2 (c e_1+ d e_2)\\
&=(-c\mu^2-2d\mu+c)e_1+(-d\mu^2 +2 c\mu+d) e_2.
\end{split}
\end{equation*}
Since the discriminant of the polynomial $-d\mu^2+2c\mu+d$ in $\mu$ is $4c^2+4d^2\geq 0$, we may find $\mu\in \bR$ such that $\hrho(e_3',e_3')\in \spa{e_1}$ and then may normalise so that $\hrho(e_3',e_3')=e_1$.

In all cases, there exist a basis $\{e_1,e_2,e_3,e_4\}$ of $\mfh$ and $\delta\in \{0,1\}$ such that $e_2=Je_1$, $e_4=Je_3$, $\hrho(e_1,\cdot)=\id_{\mfh}$, $\hrho(e_2,\cdot)=J|_{\mfh}$ and
\begin{equation*}
\hrho(e_4,e_4)=-\hrho(e_3,e_3)=-\delta e_1,\quad \hrho(e_3,e_4)=\hrho(e_4,e_3)=\delta e_2.
\end{equation*}
Thus, denoting by $\{f_1,\ldots,f_8\}$ the basis $\{e_1,\ldots,e_4,e^1,-e^2,e^3,-e^4\}$, we have
\begin{equation*}
J=\sum_{i=1}^4 \left(f^{2i-1}\otimes f_{2i}-f^{2i}\otimes f_{2i-1}\right),\quad
\omega=f^{15}-f^{26}+f^{37}-f^{48}
\end{equation*}
and in Salamon's notation, the differentials of the elements of the dual basis $\{f^1,\ldots,f^8\}$ are
\begin{equation*}
(0^4,-15+26-37+48,-16-25-38-47,-17+28-\delta(35-46),-18-27-\delta(36+45)).
\end{equation*}
\end{proof}
\subsection{Complex symplectic manifolds with Lagrangian fibrations}\label{subsec:Lag_fib}

We show how to use the complex symplectic cotangent extension to produce examples of compact complex symplectic manifolds endowed with a Lagrangian fibration. If we start with a nilpotent Lie algebra $\frh$ of dimension $2n$ endowed with a complex structure $J$, we obtain a complex symplectic structure $(\bJ,\bO)$ on the Lie algebra $\frg=\frh^*\oplus\frh$. If $\rho$ can be chosen in such a way that $\frg$ is nilpotent, then one can use Mal'tsev theorem to check whether $G$, the unique connected, simply connected, nilpotent Lie group with Lie algebra $\frg$, admits a lattice $\Gamma$. If this is the case, then $N=\Gamma\backslash G$ is a compact (nil)manifold endowed with an invariant complex symplectic structure. The short exact sequence $0\to\frh^*\to\frg\to\frh\to 0$ of Lie algebras produces a short exact sequence $1\to\bR^{2n}\to G\stackrel{\pi}{\to} H\to 1$ of Lie groups. It is known that then $\bR^{2n}\cap \Gamma\subset\bR^{2n}$ and $\pi(\Gamma)=(\Gamma\cap\bR^{2n})\backslash \Gamma\subset H$ are lattices. We obtain therefore a principal $T^{2n}$-bundle
\[
\xymatrix{
T^{2n}\ar[r] & N\ar[d]\\
& B
}
\]
whose base $B\coloneqq \pi(\Gamma)\backslash H$ is a $2n$-dimensional nilmanifold and whose fibers are Lagrangian tori. This can be compared with the structure of a Lagrangian fibration on a projective/K\"ahler irreducible holomorphic symplectic manifold: it is conjectured that the base of such a fibration is isomorphic to $\bC P^n$. The conjecture is true, even under the K\"ahler hypothesis, if the base of the Lagrangian fibration is smooth (see \cite{greb-lehn}).

\begin{example}
Consider the third set of solutions for $\rho$ in Example \ref{ex:cot_extension} above, with all parameters equal to 0. Consider 
the basis $\{f_1,\ldots,f_{12}\}=\{e_1,\ldots,e_6,e^1,\ldots e^6\}$ of $\frg=\frh^*\oplus\frh$; with respect to the dual basis $\{f^1,\ldots,f^{12}\}$ of $\frg^*$, the non-zero differential are
\[
df^4=f^{12}\,, \quad df^5=f^{13}\,, \quad df^6=f^{23}\,, \quad df^7=f^{2,10}+f^{3,11} \quad \mathrm{and} \quad df^8=f^{2,9}+f^{3,12}\,,
\]
and $\frg$ is 2-step nilpotent. The symplectic form on $\frg$ is $\bO=\sum_{i=1}^6f^{i+6}\wedge f^i$, and the complex structure is $\bJ=J^*\oplus J$. Clearly $\frg$ has a rational structure constants, hence $G$ has a lattice $\Gamma$. We thus obtain a compact, complex symplectic (nil)manifold $N=\Gamma\backslash G$ with a Lagrangian fibration
\[
\xymatrix{
T^6\ar[r] & N\ar[d]\\
& B
}
\]
where $B$ is a nilmanifold associated the Lie algebra $\frh_7$. A basis for holomorphic 1-forms is given by
\[
\omega^1=f^1-if^2\,, \ \omega^2=f^3+if^4\,, \ \omega^3=f^5-if^6\,, \ \omega^4=f^7+if^8\,, \ \omega^5=f^9-if^{10}\,, \ \omega^6=f^{11}+if^{12}\,,
\]
with non-zero differentials
\[
d\omega^2=\frac{1}{2}\omega^{1\bar 1}\,, \quad d\omega^3=\frac{1}{2}\left(\omega^{12}+\omega^{1\bar 2}\right) \quad \mathrm{and} \quad d\omega^4=\frac{1}{2}\left(-\omega^{15}+\omega^{\bar1 5}+\omega^{26}+\omega^{\bar 2 6}\right)\,.
\]
The complex symplectic form is $\bO_\bC=-(\omega^{14}+\omega^{25}+\omega^{36})$. Also,
\[
H^{2,0}(N)=\spa{\omega^{12},\omega^{13},\omega^{15},\omega^{16},\omega^{56},\bO_\bC}\,.
\]
\end{example}

\section{Abelian complex symplectic structures on solvable Lie algebras with non trivial center}\label{sec:Abelian}

In this section we investigate complex symplectic structures on solvable Lie algebras with non-trivial center (this includes for instance nilpotent Lie algebras) under the assumption that the complex structure is Abelian. Let $\mathfrak{g}$ be a Lie algebra admitting a complex symplectic structure $(J,\omega)$; assume that $J$ is Abelian. Set $\mathfrak{g}^1:=[\mathfrak{g},\mathfrak{g}]$ and $\mathfrak{g}^1_J:=\mathfrak{g}^1+J\mathfrak{g}^1$. In particular $\mathfrak{g}^1_J$ is a $J$-invariant Lie subalgebra of $\mathfrak{g}$. Recall that since $J$ is Abelian, $\mathfrak{g}^1$ is an Abelian ideal in $\mathfrak{g}$, so that $\mathfrak{g}$ is $2$-step solvable by \cite[Lemma 1]{petravchuk}.


\begin{proposition}\label{prop:J:Abelian}
Let $(\frg,J,\omega)$ be a complex symplectic Lie algebra with $J$ Abelian, and let $\mfz$ be its center. 
\begin{enumerate}[1.]
\item $(\mathfrak{g}^1)^{\perp_\omega}$ and $(\mathfrak{g}^1_J)^{\perp_\omega}$ are Abelian.
\item $(\mathfrak{g}^1_J)^{\perp_\omega}$ is $J$-invariant.
\item $\mfz$ is $J$-invariant and contained in $(\mathfrak{g}^1_J)^{\perp_\omega}$.
\item $J[X,Y]=[X,JY]$, $\forall~X\in(\mathfrak{g}^1_J)^{\perp_\omega}$.
\end{enumerate}
\end{proposition}

\begin{proof}


It suffices to show that $(\mathfrak{g}^1)^{\perp_\omega}$ is Abelian, since $\mathfrak{g}^1_J \supseteq \mathfrak{g}^1$ implies $(\mathfrak{g}^1_J)^{\perp_\omega}\subseteq (\mathfrak{g}^1)^{\perp_\omega}$. If $X,Y\in (\mathfrak{g}^1)^{\perp_\omega}$, then
\[
0=\omega([X,Y],Z)+\omega([Y,Z],X)+\omega([Z,X],Y)=\omega([X,Y],Z)
\]
for every $Z\in\mathfrak{g}$, since $\omega$ is closed; since it is also non-degenerate, $[X,Y]=0$ (notice that this holds without any assumption on $J$). This proves 1.

It is immediate to see that $(\mathfrak{g}^1_J)^{\perp_\omega}$ is $J$-invariant, since $J$ is symmetric and $\mathfrak{g}^1_J$ is $J$ invariant. Indeed, if $X\in(\mathfrak{g}^1_J)^{\perp_\omega}$ then, for every $Y\in\mathfrak{g}^1_J$
\[
\omega(JX,Y)=\omega(X,JY)=0\,.
\]
Hence we have 2.

Notice that $\mfz$ is $J$-invariant, since $J$ is Abelian; moreover, it is contained in $(\mathfrak{g}^1_J)^{\perp_\omega}$. Indeed if $X\in\mfz$, then since $\omega$ is closed,
\[
\omega([Y,Z],X)=0\,;
\]
furthermore, since $J$ is symmetric,
\[
\omega(J[Y,Z],X)=\omega([Y,Z],JX)=0
\]
for every $Y,Z\in\frg$, and we have proved 3.

Suppose now $X\in(\mathfrak{g}^1_J)^{\perp_\omega} $; then, for every $Y,Z\in\frg$,
\begin{align*}
\omega(J[X,Y],Z)&=\omega([X,Y],JZ)=-\omega([JZ,X],Y)=
\omega([Z,JX],Y)=-\omega([JX,Y],Z)-\omega([Y,Z],JX)\\
&=\omega([X,JY],Z)\,,
\end{align*}
where we used that $\omega$ is closed and $J$ is Abelian. By non-degeneracy, we conclude $J[X,Y]=[X,JY]$ for every $X\in(\mathfrak{g}^1_J)^{\perp_\omega}$, and we have 4.
\end{proof}

\begin{corollary}
Let $(\frg,J,\omega)$ be a 2-step nilpotent complex symplectic Lie algebra, and let $\mfz$ be its center.
\begin{enumerate}[1.]
\item $\mfg^1$ is an $\omega$-isotropic ideal.
\item If $J$ is Abelian, $\mfg^1_J$ is a $J$-invariant $\omega$-isotropic ideal.
\item If $J$ is Abelian, then $\mathfrak{g}^1_J\subseteq\mfz\subseteq(\mathfrak{g}^1_J)^{\perp_\omega}$
\end{enumerate}
\end{corollary}

\begin{proof}
Since $\omega$ is closed, for every $X,Y,U,V\in\mathfrak{g}$,
\[
\omega([X,Y],[U,V])=-\omega([Y,[U,V]],X)-\omega([[U,V],X],Y)=0\,,
\]
which gives 1. Assuming $J$ is Abelian,
\begin{align*}
\omega([X,Y],J[U,V])&=-\omega([Y,J[U,V]],X)-\omega([J[U,V],X],Y)=
\omega([JY,[U,V]],X)+\omega([[U,V],JX],Y)\\
&=0\,;
\end{align*}
since $J$ is symmetric, using 1.,
\[
\omega(J[X,Y],J[U,V])=-\omega([X,Y],[U,V])=0\,;
\]
thus we proved 2. For 3., the first inclusion follows directly by $2$-step nilpotency and $J$ Abelian, the second inclusion was proved in Proposition \ref{prop:J:Abelian}.
\end{proof}

Let $(\mfg,J,\omega)$ be a complex symplectic solvable Lie algebra of dimension $4n$ with $J$ Abelian and non-trivial center $\mfz$. By Proposition \ref{prop:J:Abelian}, $\mfz$ is $J$-invariant and, hence, by \cite[Corollary 4.10]{BFLM}, $(\mathfrak{g},J,\omega)$ is the \emph{complex symplectic oxidation} of a complex symplectic Lie algebra $(\bar \mfg,\bar J,\bar \omega)$ of dimension $4n-4$.

For this, let us recall from \cite{BFLM}, when a $(\mfg,J,\omega)$ complex symplectic solvable Lie algebra of dimension $4n$ is the \emph{complex symplectic oxidation} of a $4n-4$-dimensional  $(\mathfrak{g},J,\omega)$:

First of all, we have
\begin{equation*}
\mfg=V\oplus \bar\mfg\oplus V^*
\end{equation*}
as vector spaces for a two-dimensional real vector space $V$ and the non-zero Lie brackets (up to anti-symmetry) are given by
\begin{eqnarray*}
\left[(v,0,0),(w,0,0)\right] &=&(0,\nu(v,w),\tau(v,w)), \label{eq:Liebracketong1}\\
\left[(v,0,0),(0,X,0)\right] &=&(0,f(v,X),g(v,X)), \label{eq:Liebracketong2}\\
\left[(0,X,0),(0,Y,0)\right] &=&(0,[X,Y]_{\bar{\mfg}},\beta(X,Y)), \label{eq:Liebracketong3}
\end{eqnarray*}
for $v,w\in V$ and $X,Y\in \bar{\mfg}$, where 
\begin{itemize}
	\item $\nu\colon\Lambda^2 V\to \bar{\mfg}$ $(\Leftrightarrow \nu\in\Lambda^2 V^*\otimes\bar\mfg)$,
	\item $\tau\colon\Lambda^2 V\to V^*$ $(\Leftrightarrow \tau\in\Lambda^2 V^*\otimes V^*)$,
	\item $f\colon V\otimes \bar\mfg\to\bar{\mfg}$ $(\Leftrightarrow f\in V^*\otimes\bar\mfg^*\otimes\bar\mfg\cong V^*\otimes\mathrm{End}(\bar\mfg))$,
	\item $g\colon V\otimes \bar{\mfg}\to V^*$ $(\Leftrightarrow g\in V^*\otimes\bar\mfg^*\otimes V^*)$,
	\item $\beta\colon \Lambda^2 \bar{\mfg}\to V^*$ $(\Leftrightarrow \beta\in \Lambda^2\bar\mfg^*\otimes V^*)$,
\end{itemize}
and $[\cdot,\cdot]_{\bar{\mfg}}$ is the Lie bracket of $\bar{\mfg}$. Moreover, $V$ should admit an almost complex structure $I$ such that the complex structure $J$ on $\mfg=V\oplus\bar\mfg\oplus V^*$ is given by $J=I+\bar J+I^*$ and $\omega$ on $\mfg=V\oplus\bar\mfg\oplus V^*$ is given by
\begin{equation*}
\omega((v,X,\alpha),(w,Y,\beta))=\alpha(Y)-\beta(X)+\bar\omega(X,Y)
\end{equation*}
for $v,w\in V$, $X,Y\in \bar\mfg$, $\alpha,\beta\in V^*$. Now there are more conditions on all these tensors which ensure that $(J,\omega)$ is a complex symplectic structure and which we determined in detail in \cite{BFLM}. We will recall them below but first specialise to the case that $J$ is Abelian noting that this is the case if and only if
\begin{itemize}
\item $\bar{J}$ is Abelian;
\item $f(Iv,\bar J X)=f(v,X)$, for all $v\in V$, $X\in \bar\mfg$;
\item $g(Iv,\bar J X)=g(v,X)$, for all $v\in V$, $X\in \bar\mfg$;
\item $\beta$ is of type $(1,1)$.
\end{itemize}
Now we try to simplify these relations using also the mentioned conditions from \cite{BFLM} that ensure that $(J,\omega)$ is a complex symplectic structure. Note that these conditions show that already the triple $(f,S,\tau)$ is enough to determine all other tensors from above, where $S$ is the symmetric part of $g\in V^*\otimes V^*\otimes \bar\mfg$ and such a triple $(f,S,\tau)$ satisfying all necessary conditions to ensure that $(J,\omega)$ is a complex symplectic structure is called \emph{complex symplectic oxidation data} on $(\bar\mfg,\bar J,\bar \omega)$.

For the mentioned simplification, we choose $v_1\in V\setminus \{0\}$, set $v_2\coloneqq I v_1$ and $f_j\coloneqq f(v_j)$ for $j=1,2$. Moreover, we denote by $A$ the anti-symmetric part of $g$ so that $g=S+A$, and set $S_{jk}:=S(v_j,v_k)\in \bar\mfg$ and $A_{jk}:=A(v_j,v_k)\in \bar\mfg$ for $j,k=1,2$. \\

For $j=1,2$, we also set $f_j^{\bar J}\coloneqq\frac{f_j-\bar J\circ f_j\circ \bar J}{2}$ and $f_j^{-\bar J}\coloneqq\frac{f_j+\bar J\circ f_j\circ \bar J}{2}$, hence decomposing $f_j=f_j^{\bar J}+f_j^{-\bar J}$ into its $\bar J$-invariant and $\bar J$-anti-invariant part. We now are able to prove the following
\begin{proposition}\label{pro:conditionsJAbelian}
A $4n$-dimensional complex symplectic Lie algebra $(\mfg,J,\omega)$ which is the complex symplectic oxidation of a $(4n-4)$-dimensional complex symplectic Lie algebra $(\bar\mfg,\bar J,\bar\omega)$ has Abelian $J$ if and only if
\begin{itemize}
\item $\bar J$ is Abelian;
\item $f_2^{\bar J}=-\bar J f_1^{\bar J}$, $f_2^{-\bar J}=\bar J f_1^{-\bar J}$, $f_1^{\bar J}\in \mathfrak{sp}(\bar \mfg,\bar J,\bar \omega)$;
\item $S_{12}=-\frac{S_{11}-S_{22}}{2}\circ \bar J$.
\end{itemize}
\end{proposition}

\begin{proof}
First of all, notice that the above conditions $f(Iv,\bar J X)=f(v,X)$ and $g(Iv,\bar J X)=g(v,X)$ for all $v\in V$, $X\in \bar\mfg$ are equivalent to
\begin{equation}\label{eq:conditionsJAbelianprelim}
\begin{split}
f_2\circ \bar J=f_1\,,\qquad (S_{12}-A_{12})\circ \bar J=S_{11}\,,\qquad S_{22}\circ \bar J=S_{12}+A_{12}\,.
\end{split}
\end{equation}
Now by \cite[Lemma 4.11]{BFLM}, $f_2-\bar J\circ f_1\in \mathfrak{sp}(\bar g,\bar J,\bar \omega)$. 
Since $f_2=-f_1\circ \bar J$, we have $f_1^{\bar J}\in \mathfrak{sp}(\bar g,\bar J,\bar \omega)$ and
\begin{equation*}
f_2^{\bar J}=-\bar J\circ f_1^{\bar J}\,,\qquad f_2^{-\bar J}=\bar J\circ f_1^{-\bar J}\,,\qquad f_1^{\bar J},f_2^{\bar J}\in \mathfrak{sp}(\bar \mfg,\bar J,\bar \omega)\,,
\end{equation*}
where $f_2^{\bar J}=-\bar J f_1^{\bar J}\in \mathfrak{sp}(\bar \mfg,\bar J,\bar \omega)$ follows from $f_1^{\bar J}\in \mathfrak{sp}(\bar \mfg,\bar J,\bar \omega)$. Now by \cite[Proposition 4.6]{BFLM}, $\beta=-f.\bar \omega$. Since $\bar\omega$ is of type $(2,0)+(0,2)$, $f^{\bar J}.\bar\omega$ is also of type $(2,0)+(0,2)$ while $f^{-\bar J}.\bar\omega$ is of type $(1,1)$. So $f^{\bar J}.\bar\omega$ has to vanish, and this holds true, since $f_1^{\bar J},f_2^{\bar J}\in \mathfrak{sp}(\bar g,\bar J,\bar \omega)$. Next, \cite[Proposition 4.6]{BFLM} contains a formula for $A$, which, however, has a wrong factor in front. To correct this wrong factor, first note that by \cite[Lemma 4.4]{BFLM}, we have $A=\tfrac{1}{2}\nu\hook \bar\omega$. The error lies in \cite[Lemma 4.5]{BFLM}, which  states $\nu=\bar J(\mathrm{Alt}(S_I))^{\sharp}$ but the proof of that lemma actually gives
$\nu=2 \bar J(\mathrm{Alt}(S_I))$. Hence, $A= \bar J^*(\mathrm{Alt}(S_I))$, which is equivalent to  $A_{12}=\tfrac{1}{2}(S_{11}+S_{22})\circ \bar J$. The last two equations in \eqref{eq:conditionsJAbelianprelim} are then equivalent to
\begin{equation*}
S_{12}=-\frac{S_{11}-S_{22}}{2}\circ \bar J.
\end{equation*}
\end{proof}

As a consequence we prove the following
\begin{theorem}\label{thm:nilp_Abelian}
Any nilpotent complex symplectic Lie algebra $(\mfg,J,\omega)$ with Abelian $J$ may be obtained by iterated complex symplectic oxidation from the trivial $0$-dimensional complex symplectic Lie algebra.
\end{theorem}
\begin{proof}
As already stated already above, every $4n$-dimensional complex symplectic Lie algebra $(\mfg,J,\omega)$ with Abelian $J$ is the complex symplectic oxidation of a $(4n-4)$-dimensional $(\bar\mfg,\bar J,\bar\omega)$ and by Proposition \ref{pro:conditionsJAbelian}, $\bar J$ is again Abelian and so also $(\bar\mfg,\bar J,\bar\omega)$ is the complex symplectic oxidation of a complex symplectic Lie algebra of dimension $4n-8$. Iterating this process we obtain the desired result.
\end{proof}

Finally, we consider the possible step lengths of a $4n$-dimensional nilpotent complex symplectic Lie algebras $(\mfg,J,\omega)$ with Abelian $J$. For this, we first note that then $J$ is also nilpotent. Consequently, the maximal nilpotency step of $\mfg$ is $2n$. We show now that any nilpotency step is possible. Even more, we prove:
\begin{proposition}\label{pro:nilpotencystepAbelian}
Let $n\in \bN$ be fixed. Then, for any $m\in \{1,\ldots,2n\}$, there exists a $4n$-dimensional nilpotent complex symplectic Lie algebra $(\mfg,J,\omega)$ of step length $m$. Moreover, if $m\geq 2$ and $(n,m)\neq (1,2)$ we may find $(\mfg,J,\omega)$ as above with $\mfg$ being indecomposable and both with $J$ being non-Abelian nilpotent or with $J$ being Abelian.
\end{proposition}
\begin{remark}
	Consider the case $(n,m)=(1,2)$, i.e. a $4$-dimensional nilpotent complex symplectic Lie algebra $(\mfg,g,J)$ of step length $2$. By
	 \cite[Proposition 5.7]{BFLM}, any such complex symplectic Lie algebra $(\mfg,g,J)$ is isomorphic to  $(\mfh_3\oplus\bR,J_0,\omega_0)$
	with a ``standard'' complex symplectic structure $(J_0,\omega_0)$ on $\mfh_3\oplus \bR$ with Abelian $J_0$. Hence, for any $4$-dimensional nilpotent complex symplectic Lie algebra $(\mfg,J,\omega)$ of step length $2$, $\mfg$ is decomposable
	and $J$ is Abelian.
\end{remark}
\begin{proof}[Proof of Proposition \ref{pro:nilpotencystepAbelian}]
We first assume that $m\geq 3$ and treat the cases $m=2l+1$ or $m=2l+2$ for some $l\in \{1,\ldots,2n-1\}$ simultaneously:

In this case, we obtain the desired examples by complex symplectic oxidation of the Abelian complex symplectic Lie algebra
$(\bR^{4n-4},J_0,\omega_0)$ with
\begin{equation*}
\begin{split}
J_0&:=\sum_{j=1}^{2(n-1)} -e^{2j-1}\otimes e_{2j}+e^{2j}\otimes e_{2j-1},\\
\omega_0&:=\sum_{j=1}^{l} (-1)^{j-1} \left(e^{2j-1}\wedge e^{2(2l-j+1)}+e^{2j}\wedge e^{2(2l-j+1)-1}\right)+\sum_{k=l+1}^{n-1} \left(e^{4k-3}\wedge e^{4k}+e^{4k-2}\wedge e^{4k-1}\right)
\end{split}
\end{equation*}
In all cases, we set $\tau:=0$, choose a specific $f_1\in \mathfrak{sp}(\bR^{4n-4},J_0,\omega_0)$ and then set $f_2\coloneqq-J_0 f_1\in \mathfrak{sp}(\bR^{4n-4},J_0,\omega_0)$. Moreover, recall that \cite[Proposition 4.6]{BFLM} describes the necessary conditions for the above tensors on $\bar \mfg$ to give rise to a complex symplectic Lie algebra $(\mfg,J,\omega)$ in four more dimensions and that conditions 1. 6. and 7. in this proposition are satisfied independently of our choice of $S$ since by our choices $f_1,f_2$ commute and preserve $\omega_0$. Moreover, condition 2. in \cite[Proposition 4.6]{BFLM} is given by $\beta=-f.\omega_0=0$.

Explicitly, we set
\begin{equation*}
f_1(e_{2j-1}):=e_{2j+1},\qquad f_2(e_{2j}):=e_{2j+2}
\end{equation*}
for all $j=1,\ldots,2l-1$
\begin{equation*}
f_1(e_{4l-1}):=f_1(e_{4l}):=0
\end{equation*}
and
\begin{equation*}
f_1(e_{4k-3}):=e_{4k-1},\qquad f_1(e_{4k-2}):=e_{4k},\qquad f_1(e_{4k-1}):=f_1(e_{4k}):=0
\end{equation*}
for all $k=l+1,\ldots,n$. Then
\begin{equation*}
\begin{split}
f_2(e_{2j-1})&=e_{2j+2},\qquad f_2(e_{2j})=-e_{2j+1},\qquad f_2(e_{4l-1})=f_2(e_{4l})=0,\\
f_2(e_{4k-3})&=e_{4k},\qquad f_2(e_{4k-2})=-e_{4k-1},\qquad f_2(e_{4k-1})=f_2(e_{4k-2})=0
\end{split}
\end{equation*}
for all $j=1,\ldots,2l-1$ and all $k=l+1,\ldots,n$.

Now we distinguish the cases $m=2l+1$ and $m=2l+2$:
\begin{itemize}
	\item $m=2l+1$:
	Set $S_{11}:=-S_{22}:=e^{4l-1}$. By conditions 3. and 4. in \cite[Proposition 4.6]{BFLM}
	we have
	\begin{equation*}
	\nu(v_1,v_2)=J_0 (S_{11}+S_{22})^{\sharp}=0,A_{12}=\tfrac{1}{2}(S_{11}+S_{22})\circ J_0=0.
	\end{equation*}
	Now condition 5. in \cite[Proposition 4.6]{BFLM} equals $\mathrm{Alt}(S\circ f)=0$, which is equivalent to
	\begin{equation*}
	S_{12}\circ f_1=S_{11}\circ f_2=-e^{4l-2},\quad S_{12}\circ f_2=S_{22}\circ f_1=-e^{4l-3}.
	\end{equation*}
	These equations are fulfilled for $S_{12}=-e^{4l}$ and also for $S_{12}=-e^{4l}+e^1$.
	In the first case, we also have $S_{12}=-\tfrac{1}{2}(S_{11}-S_{22})\circ J_0$ whereas in the second case, this equation does not hold. Hence, Proposition \ref{pro:conditionsJAbelian} yields that in the first case the corresponding complex symplectic oxidation $(\mfg,J,\omega)$ has Abelian $J$, whereas in the second case, $J$ is not Abelian. In both cases, $\mfg$ is indecomposable and we have
	\begin{equation*}
	\begin{split}
	\mfg_1&=V^*\oplus \spa{e_{4k-1},e_{4k}|k=l+1,\ldots,2n},\\
	\mfg_2&=V^*\oplus \spa{e_{4l-1},\ldots,e_{4n}},\\
	\mfg_r&=V^*\oplus \spa{e_{2(2l-r+2)-1},\ldots,e_{4n}}, r=3,\ldots,2l-1,\\
	\mfg_{2l}&=V^*\oplus \spa{e_3,\ldots,e_{4n}},\qquad \mfg_{2l+1}=\mfg.
	\end{split}
	\end{equation*}
	This shows that $\mfg$ has step length $m=2l+1$ and  $J$ is nilpotent in both cases.
	\item $m=2l+2$:
    Here, we set $S_{11}:=S_{22}:=e^{4l-1}$. Conditions 3. and 4. in \cite[Proposition 4.6]{BFLM} then yield
	\begin{equation*}
	\nu(v_1,v_2)=J_0 (S_{11}+S_{22})^{\sharp}=2 e_1,\qquad A_{12}=\tfrac{1}{2}(S_{11}+S_{22})\circ J_0=e^{4l}.
	\end{equation*}
	We now look at condition 5. in \cite[Proposition 4.6]{BFLM} and note that
	\begin{equation*}
	\begin{split}
&S_{11}\circ f_2+\omega_0(f_1(J_0(S_{11}+S_{22})^{\sharp}),\cdot)+\tfrac{3}{2} (S_{11}+S_{22})\circ J_0\circ f_1\\
=\ &-e^{4l-2}+\omega_0(f_1(2e_1),\cdot)+3 e^{4l-1}\circ J_0\circ f_1=-e^{4l-2}+2\omega_0(e_3,\cdot)+3 e^{4l-2}\\
=\ &-e^{4l-2}-2e^{4l-2}+3 e^{4l-2}=0
\end{split}
\end{equation*}
and
\begin{equation*}
\begin{split}
&-S_{22}\circ f_1+\omega_0(f_2(J_0(S_{11}+S_{22})^{\sharp}),\cdot)+\tfrac{3}{2} (S_{11}+S_{22})\circ J_0\circ f_2\\
=\ &-e^{4l-3}+\omega_0(f_2(2e_1)),\cdot)+3 e^{4l-1}\circ J_0\circ f_2=-e^{4l-3}+2\omega_0(e_4,\cdot)+3 e^{4l-3}\\
=\ &-e^{4l-3}-2 e^{4l-3}+3 e^{4l-3}=0.
\end{split}
\end{equation*}
Hence, condition 5. in \cite[Proposition 4.6]{BFLM} is satisfied if and only if
\begin{equation*}
S_{12}\circ f_1=S_{12}\circ f_2=0,
\end{equation*}
i.e. if $S_{12}\in \spa{e^1,e^2}$. We consider the cases $S_{12}=0$ or $S_{12}=e^1$. In the first case,
we have $S_{12}=-\tfrac{1}{2}(S_{11}-S_{22})\circ J_0$ whereas in the second case, this equation does not hold. Thus, again by Proposition \ref{pro:conditionsJAbelian}, in the first case the corresponding complex symplectic oxidation $(\mfg,J,\omega)$ has Abelian $J$, whereas in the second case, $J$ is not Abelian. In both cases, $\mfg$ is indecomposable and we have
	\begin{equation*}
	\begin{split}
	\mfg_1&=V^*\oplus \spa{e_{4k-1},e_{4k}|k=l+1,\ldots,2n},\\
	\mfg_2&=V^*\oplus \spa{e_{4l-1},\ldots,e_{4n}},\\
	\mfg_r&=V^*\oplus \spa{e_{2(2l-r+2)-1},\ldots,e_{4n}}, r=3,\ldots,2l-1,\\
	\mfg_{2l}&=V^*\oplus \spa{e_3,\ldots,e_{4n}},\\
	\mfg_{2l+1}&=V^*\oplus \bR^{4n-4},\qquad \mfg_{2l+2}=\mfg.
	\end{split}
	\end{equation*}
	Hence, in both cases, $\mfg$ has step length $m=2l+2$ and $J$ is nilpotent.
\end{itemize}
Now let us look at the case $m=2$ and so $n\geq 2$.

We first construct an example with non-Abelian nilpotent $J$. To do this, take the complex symplectic Lie algebra $(\bR^{4(n-1)},J_0,\omega_0)$ with $J_0$ as above and $\omega_0$ as in the case $l=0$. Also define, $f_1$
as in the case $l=0$ but set $f_2:=0$. Moreover, $S_{11}:=-S_{22}:=-e^2$ and $S_{12}:=e^1$ and $\tau:=0$. One checks, using \cite[Proposition 4.6]{BFLM}, that this defines complex symplectic oxidation data and that $\beta=0$, $\nu=0$ and $A=0$, and so the non-zero Lie brackets (up to anti-symmetry) on $V\oplus \spa{e_1,e_2,e_3,e_4}\oplus V^*$ are given by
\begin{equation*}
[v_1,e_1]=e_3+v^2,\quad [v_1,e_2]=e_4-v^1,\quad [v_2,e_1]=v^1,\quad [v_2,e_2]=v^2.
\end{equation*}
Hence, $\mfg$ is indecomposable as a Lie algebra with
\begin{equation*}
\mfg_1=V^*\oplus \spa{e_{4k-1},e_{4k}|k=1,\ldots,n},\qquad \mfg_2=\mfg,
\end{equation*}
i.e. $\mfg$ has step length $2$. Moreover, by Proposition \ref{pro:conditionsJAbelian}, the obtained complex symplectic Lie algebra $(\mfg,J,\omega)$ has non-Abelian $J$ due to $S_{12}=e^1\neq -e^1=-\tfrac{1}{2}(S_{11}-S_{22})\circ J_0$.

Next, we construct an example with Abelian $J$. For this, we take the complex symplectic Lie algebra $(\mfh_3\oplus \bR^{4n-7},J_0,\omega_0)$ with basis $e_1,\ldots,e_{4n}$ of $\mfh_3\oplus \bR^{4n-3}$
with the only non-zero Lie bracket (up to anti-symmetry) given by $[e_1,e_2]=e_3$. Moreover,
we define $J_0$ as in the case $\bR^{4n-4}$ and $\omega_0$ as in the case $\bR^{4n-4}$ and $l=0$. Next, we set $f_1$ and $f_2=-J_0f_1$ as in the case $\bR^{4n-4}$ and $l=0$ noting that $f_1$ and $f_2$ are both derivations of $\mfh_3\oplus \bR^{4n-7}$. Moreover, we set $\tau(v_1,v_2)\coloneqq v^1$, $S_{11}\coloneqq -S_{22}\coloneqq e^2$ and $S_{12}\coloneqq e^1$. One checks that this defines complex symplectic oxidation data on $(\mfh_3\oplus \bR^{4n-7},J_0,\omega_0)$. Since
$S_{12}=-\tfrac{1}{2}(S_{11}-S_{22})\circ J_0$ holds, Proposition \ref{pro:conditionsJAbelian} tells us that the obtained complex symplectic Lie algebra $(\mfg,J,\omega)$ has Abelian $J$. Now the only non-zero Lie brackets (up to anti-symmetry) on $V\oplus \spa{e_1,e_2,e_3,e_4}\oplus V^*$ are
\begin{equation*}
\begin{split}
[v_1,e_1]&=-[v_2,e_2]=e_3+v^2,\qquad [v_1,e_2]=[v_2,e_1]=e_4+v^1,\\
[v_1,v_2]&=v^1,\qquad [e_1,e_2]=e_3.
\end{split}
\end{equation*}
Hence, $\mfg$ is indecomposable as a Lie algebra and
\begin{equation*}
\mfg_1=V^*\oplus \spa{e_{4k-1},e_{4k}|k=1,\ldots,n},\qquad \mfg_2=\mfg,
\end{equation*}
i.e. $\mfg$ is of step length $2$.
\end{proof}


\end{document}